\documentclass[12pt]{article}
\usepackage{amsmath}
\usepackage{amssymb}
\usepackage{amscd}
\newtheorem{lemme}{Lemme}
\newtheorem{theoreme}{Th\'eor\`eme}
\newtheorem{proposition}{Proposition}
\newtheorem{corollaire}{Corollaire}
\newcounter{numeroexemple}
\newcounter{numeroremarque}
\newenvironment{exemple}
  {\addtocounter{numeroexemple}{1} 
  \begin{trivlist}\item[]\textbf{Exemple \thenumeroexemple}}{\end{trivlist}}
\newenvironment{remarque}
  {\addtocounter{numeroremarque}{1} 
\begin{trivlist}\item[]\textbf{Remarque \thenumeroremarque}}{\end{trivlist}}
\newcounter{numerodefinition}

\newcounter{numeroquestion}

\newcounter{numeroconjecture}

\newenvironment{preuve}{\begin{trivlist}\item[]\textit{Preuve.}}
{\item[] \hfill $\square$\end{trivlist}}
\def\mult{{\mbox{\rm mult}}}
\def\m{{\mbox{\rm m}}}
\def\Id{{\mbox{\rm Id}}}
\def\O{{\mbox{\rm O}}}
\def\o{{\mbox{\rm o}}}
\def\hfld#1#2{\smash{\mathop{\hbox to 12mm{\rightarrowfill}}
     \limits^{\scriptstyle#1}_{\scriptstyle#2}}}
\def\hflg#1#2{\smash{\mathop{\hbox to 12mm{\leftarrowfill}}
     \limits^{\scriptstyle#1}_{\scriptstyle#2}}}

\begin{document}
\title{Sur les endomorphismes polynomiaux permutables de
  $\mathbb{C}^2$
\footnote{{\bf Classification math\'ematique:} 30D05, 58F23.
\ \ \ \ \ \ \ \ \ \ \ \ \ \ \ \ \ \ \ \ \ \ \ \ \ \ \
\break
{\bf Mots cl\'es:} it\'er\'es, applications permutables,
  orbifold, critiquement finie.}
}
\author{Tien-Cuong Dinh
\footnote{
Math\'ematique-B\^atiment 425,
Universit\'e Paris-Sud, 91405 ORSAY Cedex (France).
\break
E-mails: TienCuong.Dinh@math.u-psud.fr.
}
}
\maketitle
%
%\begin{resume}
Dans cet article, nous d\'eterminons tous les couples 
  d'endomorphismes polynomiaux permutables de degr\'es sup\'erieurs \`a
  1 de $\mathbb{C}^2$ qui se
  prolongent en des endomorphismes holomorphes de $\mathbb{P}^2$
  et qui poss\`edent deux suites d'it\'er\'es disjointes.
%\end{resume}
%
%\english
\begin{abstract} 
\begin{center}{\bf On the commuting polynomial endomorphisms of
  $\mathbb{C}^2$}
\end{center}
We determine all couples of commuting polynomial 
endomorphisms of $\mathbb{C}^2$ that extends to holomorphic
endomorphisms of $\mathbb{P}^2$ and that have disjoint sequences of iterates.
\end{abstract}
\section{Introduction}
Soient $X$ une vari\'et\'e complexe et $f_1$, $f_2$ deux
endomorphismes holomorphes de $X$. On consid\`ere l'\'equation
fonctionnelle suivante:
\begin{equation}
f_1\circ f_2=f_2\circ f_1
\end{equation}
Cette \'equation poss\`ede des solutions triviales
$f_i=h^{n_i}:=h\circ h\circ \cdots \circ h$ ($n_i$ fois) o\`u $h$ est
un endomorphisme de $X$. Par la
suite, on suppose que:
\begin{equation}
f_1^n\not =f_2^m\mbox{ pour tout } (n,m)\not =(0,0) 
\end{equation}
Soit $\sigma$ un automorphisme holomorphe de $X$. Le couple 
$(f_1,f_2)$ v\'erifie
(1) si et seulement si $(\sigma^{-1} \circ f_1
\circ\sigma,\sigma^{-1}\circ f_2\circ\sigma)$ la v\'erifie. On dit que
le couple $(\sigma^{-1} \circ f_1
\circ\sigma,\sigma^{-1}\circ f_2\circ\sigma)$ est {\it conjugu\'e} au
couple $(f_1,f_2)$.
\par
Dans $\mathbb{C}$, si les $f_i$ sont des polyn\^omes de degr\'es $d_i\geq
2$, le probl\`eme (1)(2) est r\'esolu par Julia et Fatou
\cite{Julia,Fatou}. Pour
$f_i$ les fractions rationnelles de degr\'es
$d_i\geq 2$ de $X=\mathbb{P}^1$, ce probl\`eme est r\'esolu par Ritt et par
Eremenko \cite{Ritt,Eremenko} ({\it voir} \'egalement \cite{LevinPrzytycki}). 
Notons $w=[w_0:w_1:\cdots:w_k]$ les coordonn\'ees homog\`enes
de $\mathbb{P}^k$ et $z=(z_1,z_2,\ldots, z_k)$ ses coordonn\'ees
afines o\`u $z_i:=w_i/w_0$. 
Voici les solutions de (1)(2) pour $X=\mathbb{P}^1$ 
\`a une classe de conjugaison pr\`es:
\begin{enumerate}
\item $(f_1,f_2)=(z^{\pm d_1}, \lambda z^{\pm d_2})$ pour une certaine 
  racine de l'unit\'e $\lambda$;
\item $(f_1,f_2)=(\pm T_{d_1},\pm T_{d_2})$ o\`u $T_d$ est le polyn\^ome
  de Tchebychev de degr\'e $d$ d\'efini par $T_d(\cos t):=\cos dt$. Les
  signes $\pm$ sont choisis selon la parit\'e des $d_i$. 
Les solutions sont les quatres couples 
$(\pm T_{d_1},\pm T_{d_2})$ pour les $d_i$ impairs,
  un seul couple $(T_{d_1},T_{d_2})$ (\`a une classe de conjugaison
  pr\`es) pour les autres cas;
\item Il existe un tore complexe
  $\mathbb{T}=\mathbb{C}/\Gamma$, une fonction elliptique 
  $F: \mathbb{T}\rightarrow \mathbb{P}^1$, des applications lin\'eaires
  $\Lambda_i(z)=a_iz+b_i$ qui pr\'eservent le groupe
  $\Gamma$ et v\'erifient $F\circ f_i=F\circ \Lambda_i$,
  $\Lambda_1\circ \Lambda_2=\Lambda_2\circ \Lambda_1$ sur $\mathbb{T}$.
La fonction $F$ et les applications $\Lambda_i$ sont d\'ecrites dans
  \cite{Eremenko}. 
\end{enumerate} 
Nous donnons ci-dessous quelques exemples pour
$X=\mathbb{C}^2$. Pour simplifier les notations, $T_d z_i$
signifiera la valeur de $T_d$ en $z_i$. 
\begin{exemple} $f_1(z_1,z_2):=(z_1^{d_1},
  \pm T_{d_1}z_2)$ et $f_2(z_1,z_2):=(\lambda z_1^{d_2},
  \pm T_{d_2}z_2)$ 
sont permutables, o\`u $\lambda$ est une certaine racine de l'unit\'e
 et les signes $\pm$ sont choisis selon la parit\'e des
  $d_i$ ({\it voir} le cas $X=\mathbb{P}^1$).
\end{exemple}
\begin{exemple} $f_i(z_1,z_2):=(\pm
  T_{d_i}z_1, \pm T_{d_i} z_2)$ ou $(\pm T_{d_i}z_2, \pm T_{d_i}z_1)$
  sont permutables o\`u les signes $\pm$ sont choisis selon 
  la parit\'e des $d_i$.
\end{exemple}
\begin{exemple} Soient $P_i$, $Q_i$ les polyn\^omes homog\`enes \`a
  deux variables, de
  degr\'e $d_i$ et sans facteur commun. Supposons que les  $[P_i:Q_i]$
  d\'efinissent sur $\mathbb{P}^1$ deux endomorphismes
  permutables. Alors il existe des constantes $\lambda_i\not = 0$ 
  telles que les $f_i:=(\lambda_i P_i,\lambda_iQ_i)$ d\'efinissent des
  endomorphismes holomorphes permutables de $\mathbb{C}^2$ qui se
  prolongent holomorphiquement \`a l'infini en des endomorphismes
  permutables de $\mathbb{P}^2$. 
\end{exemple}
\begin{exemple} Soient $h_1$, $h_2$ deux endomorphismes holomorphes
  de $\mathbb{P}^1$ et $\pi$ l'application holomorphe de
  $\mathbb{P}^1\times \mathbb{P}^1$ dans $\mathbb{P}^2$ qui d\'efinit
  un rev\^etement ramifi\'e \`a deux feuillets et v\'erifie
  $\pi(x,y)=\pi(y,x)$. Si $([x_0:x_1],[y_0:y_1])$ sont les
  coordonn\'ees de $\mathbb{P}^1\times \mathbb{P}^1$, on a
  $\pi([x_0:x_1],[y_0:y_1])=[x_0y_0:x_0y_1+x_1y_0:x_1y_1]$.
Soient $F_i(a,b):=(h_ia,h_ib)$ deux
  endomorphismes de $\mathbb{P}^1\times\mathbb{P}^1$. Il existe des
  endomorphismes holomorphes $f_i$ de $\mathbb{P}^2$ tels que
  $f_i\circ\pi=\pi\circ F_i$. Lorsque $h_1,h_2$ sont
  permutables, les $f_i$ sont permutables. Lorsque les
  $h_i$ sont des polyn\^omes, 
  $f_1$, $f_2$ sont polynomiaux.
\end{exemple}
De m\^eme mani\`ere, on peut construire des
endomorphismes permutables de $\mathbb{P}^k$. Pour $k\geq 3$, 
d'autres exemples sont construits dans
\cite{Veselov}. Les automorphismes polynomiaux permutables de
$\mathbb{C}^2$ sont \'etudi\'es dans \cite{Lamy}.  
\par
Notons ${\cal PH}_d$ l'ensemble des endomorphismes polynomiaux de
$\mathbb{C}^2$ qui se prolongent holomorphiquement \`a l'infini. 
Notre r\'esultat principal est:
\begin{theoreme} Soient $f_1\in {\cal PH}_{d_1}$ et $f_2\in {\cal PH}_{d_2}$ 
v\'erifiant les conditions (1), (2)
o\`u $d_1\geq 2$ et $d_2\geq 2$.
  Alors le couple $(f_1,f_2)$ est conjugu\'e \`a l'un des couples
  d\'ecrits dans les exemples 1-4.
\end{theoreme}
Un {\it orbifold} est un couple ${\cal O}=(X,n)$, o\`u $X$ est une vari\'et\'e
complexe et $n$ est une fonction \`a valeur dans
$\mathbb{N}_+\cup\{\infty\}$, d\'efinie sur
l'ensemble des hypersurfaces irr\'eductibles de $X$; cette fonction
est \'egale \`a 1 sauf sur un ensemble localement
fini d'hypersurfaces. {\it Une
application holomorphe} d'un orbifold ${\cal O}$ dans un autre ${\cal
  O}'=(X',n')$ de m\^eme dimension 
est une application holomorphe ouverte $f$ de $X\setminus H$ dans
$X'$ telle que $n'(f({\cal D}))=\mult(f,{\cal D})n({\cal D})$ 
pour toute hypersurface
irr\'eductible ${\cal D}$ de $X$, o\`u $H$ est la r\'eunion des hypersurfaces
$H_j$ v\'erifiant $n(H_j)=\infty$ et $\mult(f,{\cal D})$ 
est la multiplicit\'e de
$f$ en un point g\'en\'erique de ${\cal D}$ (on dira simplement {\it la
multiplicit\'e} de $f$ sur ${\cal D}$). 
Lorsque $f$ est un rev\^etement ramifi\'e au-dessus de
$X'\setminus H'$, on dit que $f$ d\'efinit
{\it un rev\^etement} ${\cal O}$
au-dessus de ${\cal O}'$, o\`u $H'$ est la r\'eunion des hypersurfaces
$H'_j$ v\'erifiant $n'(H'_j)=\infty$. Remarquons que si $f$ d\'efinit
un rev\^etement d'un orbifold dans lui-m\^eme, alors l'ensemble
critique ${\cal C}$ de $f$ est pr\'ep\'eriodique, i.e. $f^{n}{\cal
  C}=f^m{\cal C}$ pour certains $0\leq n<m$. Une telle application
s'appelle {\it critiquement finie}. Si $f^n$ d\'efinit un rev\^etement de
${\cal O}$ dans lui-m\^eme, $f$ d\'efinit \'egalement un rev\^etement
de ${\cal O}$ dans lui-m\^eme.
\par
Dans \cite{Eremenko}, Eremenko a montr\'e que si $f_i$ v\'erifient
(1)(2) pour $X=\mathbb{P}^1$, il existe un orbifold ${\cal
  O}=(\mathbb{P}^1,n)$ 
tel que les $f_i$ d\'efinissent des rev\^etements 
de ${\cal O}$ dans lui-m\^eme. On d\'eduit de ce
r\'esultat et du th\'eor\`eme 1 que:  
\begin{corollaire} Sous l'hypoth\`ese du th\'eor\`eme 1, il existe un
  orbifold $\mathcal{O}=(\mathbb{P}^2,n)$ tel que les applications
  $f_i$ d\'efinissent des rev\^etements de $\mathcal{O}$ dans lui-m\^eme.
En particulier, les applications $f_i$ sont critiquement finies.
\end{corollaire}
Dans le paragraphe 2, nous donnons quelques propri\'et\'es des
suites d'it\'er\'es de $f_1$ et $f_2$. 
On montre que
si $(f_1^m,f_2^n)$ appartient aux exemples 1--4, $(f_1,f_2)$ sera
conjugu\'e \`a l'un des couples d\'ecrits 
dans ces exemples (lemme 2). Par cons\'equent,
au cours de la preuve du th\'eor\`eme 1, on peut toujours remplacer
les $f_i$ par leurs it\'er\'es. On note $f_i|_\infty$ la
restriction de $f_i$ \`a la droite infinie $\{w_0=0\}$ et $\tilde
T_d(x):=2T_d(x/2)$. D'apr\`es le lemme 2, 
il suffit de consid\'erer les
cas suivants: 
$(f_1|_\infty,f_2|_\infty)=(x^{d_1},x^{d_2})$,  
$(f_1|_\infty,f_2|_\infty)=(\tilde T_{d_1}, \tilde T_{d_2})$ et le cas
o\`u les $f_i^n|_\infty$ ne sont pas conjugu\'es aux polyn\^omes pour
tout $n\geq 1$.
\par
Dans le paragraphe 3, nous donnons quelques propri\'et\'es de
l'ensemble critique d'une application holomorphe ouverte et de leur
image. 
Nous donnons \'egalement le comportement asymptotique de $f_i$ et
son ensemble critique au
point super-attractif \`a l'infini lorsque ce point existe.
\par
La preuve du th\'eor\`eme 1 se constitue par 3 derniers
paragraphes correspondant au 3 cas cit\'es ci-dessus. Dans les
paragraphes 2-5 on note ${\cal C}_i\cup\{w_0=0\}$ 
l'ensemble critique de
$f_i$ et on \'ecrit en coordonn\'ees affines 
$f_i=(f_{i1},f_{i2})$. Alors ${\cal C}_i$ est de degr\'e
$2d_i-2$ compt\'e avec les multiplicit\'es.  Pour simplifier les
notations, $fa$ et $f^{-1}a$ signifieront l'image et l'image
r\'eciproque de $a$ par $f$. 
\par
Je tiens \`a remercier Nessim Sibony de m'avoir propos\'e ce
probl\`eme ainsi que pour ses indications.
\section{Quelques propri\'et\'es de suites d'it\'er\'es}
\begin{lemme} Soient $f_1\in {\cal PH}_{d_1}$ et $f_2\in{\cal
    PH}_{d_2}$.
Supposons que $f_1\circ f_2=f_2\circ f_1$ et
  $f_1^n|_\infty=f_2^m|_\infty$. Alors il existe $n', m'\geq 1$ tels que 
$f_1^{n'}=f_2^{m'}$.  
\end{lemme}
\begin{preuve} Il suffit de consid\'erer $n=m=1$. Comme
  $f_1|_\infty=f_2|_\infty$,
  $d_1=d_2$. Posons $d:=d_1=d_2$. On peut \'ecrire les $f_i$ sur les
  coordonn\'ees homog\`enes:
$$f_1(w):=[w_0^d:P(w_1,w_2)+w_0R(w):Q(w_1,w_2)+w_0S(w)]$$
$$f_2(w):=[w_0^d:\lambda P(w_1,w_2)+w_0\tilde R(w):\lambda
Q(w_1,w_2)+w_0 \tilde S(w)]$$
o\`u $P,Q$ sont homog\`enes de degr\'e $d$, $R,S$, $\tilde R,\tilde
S$ sont homog\`enes de degr\'e $d-1$ et $\lambda\not=0$. Comme
$f_1\circ f_2=f_2\circ f_1$, $\lambda=\lambda^d$.  
Quitte
\`a remplacer $f_i$ par $f_i^{d-1}$, on peut
supposer que $\lambda=1$.
Supposons que $f_1\not =f_2$. Soit $0<\alpha\leq d$ le nombre
naturel minimal tel que:
$$f_i(w)=[w_0^d:P^*+w_0^\alpha
R_i^*+w_0^{\alpha+1}U_i:Q^*+w_0^\alpha
S_i^*+w_0^{\alpha+1}T_i]$$
o\`u $P^*, Q^*$ sont homog\`enes de degr\'e $d$ et 
ne contiennent aucun mon\^ome dont la puissance de
$w_0$ est $\geq\alpha$, $R_i^*$, $S_i^*$ sont
homog\`enes de degr\'e $d-\alpha$, ind\'ependants de $w_0$ et $U_i$,
$T_i$ sont homog\`enes de degr\'e $d-\alpha-1$.
Consid\'erons deux polyn\^omes suivants en variable $w_0$ et \`a
coefficients dans 
$\mathbb{C}[w_1,w_2]$:
$$P(P^*+w_0^\alpha R_i^*, Q^*+w_0^\alpha S_i^*)-P(P^*,Q^*).$$
Comme $f_1\circ f_2=f_2\circ f_1$, les coefficients de $w_0^\alpha$
dans ces polyn\^omes sont \'egaux:
$$P_1(P,Q) R_1^*+P_2(P,Q)S_1^*= P_1(P,Q)R_2^*+P_2(P,Q)S_2^*$$
o\`u $P_1:=\partial P/\partial w_1$ et $P_2:=\partial P/\partial w_2$.
Par cons\'equent,
$$P_1(P,Q)(R_1^*-R_2^*)+P_2(P,Q)(S_1^* -S_2^*)=0.$$
De m\^eme, 
$$Q_1(P,Q)(R_1^*-R_2^*)+Q_2(P,Q)(S_1^* -S_2^*)=0.$$
Comme $\alpha$ est maximal, $R_1^*-R_2^*$ et $S_1^*-S_2^*$ 
ne sont pas tous nuls. On en d\'eduit que $P_1/P_2=Q_1/Q_2$. 
L'application $f_i|_\infty=[P:Q]$ 
est donc constante. C'est la contradiction recherch\'ee. 
\end{preuve}
Soit $f$ un endomorphisme holomorphe de degr\'e $d\geq 2$
de $\mathbb{P}^k$. 
On note $E(f)$ l'hypersurface (\'eventuellement vide ou
r\'eductible) maximale et compl\`etement invariante par $f$,
i.e. $f^{-1}E(f)=E(f)$. Cette hypersurface existe et s'appelle
{\it l'hypersurface
exceptionelle} de $f$ \cite{Sibony} et $E(f)=E(f^n)$ pour tout
$n\geq 1$. Dans $\mathbb{P}^1$,
$\#E(f)\leq 2$; si $\#E(f)=2$, $f$ est conjugu\'e \`a $z^{\pm d}$; si
$\#E(f)=1$, $f$ est conjugu\'e \`a un polyn\^ome.
Dans $\mathbb{P}^2$, $\deg E(f)\leq 3$;  
si $\deg E(f)=1$, $f$ est conjugu\'e \`a un endomorphisme polynomial;  si
$\deg E(f)=2$, $f$ est conjugu\'e \`a $[w_p^d: w_q^d: P]$ o\`u
$\{p,q\}$ est une permutation de $\{0,1\}$ et $P$ est un polyn\^ome
homog\`ene de degr\'e $d$; si
$\deg E(f)=3$, $f$ est conjugu\'e \`a $[w_p^d: w_q^d: w_r^d]$ o\`u
$\{p,q,r\}$ est une permutation de $\{0,1,2\}$ \cite{FornessSibony}.
\begin{lemme} Sous l'hypoth\`ese du th\'eor\`eme 1, si
$(f_1^m,f_2^n)$ se trouve dans les
  exemples 1-4, $(f_1,f_2)$ est conjugu\'e \`a l'un
  des couples d\'ecrits dans ces exemples.
\end{lemme}
\begin{preuve}
\underline{\it Cas 1.}
Supposons que $(f_1^m,f_2^n)$ se trouve dans l'exemple 1. Alors l'ensemble
exceptionel de $f_i$ est la r\'eunion de deux droites dont l'une est l'infini. 
Quitte \`a changer les coordonn\'ees, on peut supposer que
$f_1=(z_1^{d_1},P(z))$ et 
$f_2=(\lambda z_1^{d_2}, Q(z))$ avec $\lambda\not =0$. 
Comme les deuxi\`emes
composantes de $f_1^m, f_2^n$ sont ind\'ependantes de $z_1$, le
polyn\^omes $P, Q$
sont ind\'ependants de $z_1$. On a donc $P\circ Q=Q\circ P$ et
$P^m=\pm T_{d_1^m}$ et $Q^n=\pm T_{d_2^n}$. Alors $(P,Q)$
est conjugu\'e \`a $(\pm T_{d_1},\pm T_{d_2})$ et $(f_1,f_2)$ est
  conjugu\'e au couple d\'ecrit dans l'exemple 1.
\\
\underline{\it Cas 2.}
Supposons que $(f_1^m,f_2^n)$ se trouve dans l'exemple 2. On a
$f_1^{2m}=(T_{d_1^{2m}}z_1,T_{d_1^{2m}}z_2)$ et
$f_2^{2n}=(T_{d_2^{2n}}z_1,T_{d_2^{2n}}z_2)$. Posons $d:=d_2^{2n}$ et 
$f:=f_2^{2n}=(T_dz_1,T_d z_2)$. Il existe un point $a\in\mathbb{C}^2$ 
fixe pour $f_1$ et p\'eriodique pour
  $f$ car l'ensemble des points fixes de $f_1$ est invarant par $f$. 
Quitte \`a remplacer $f$ par $f^l$ on peut
supposer que $a$ est fixe pour $f$. Comme $f=(T_dz_1,T_dz_2)$,  il
existe $(t_1^*,t_2^*)\in \mathbb{R}^2$ tel que $a=(\cos t_1^*,\cos
t_2^*)=(\cos dt_1^*,\cos dt_2^*)$. 
Posons $\varphi(t_1,t_2):=(\cos(t_1-t_1^*),\cos(t_2-t_2^*))$ une
application holomorphe 
d'un voisinage de $0\in\mathbb{C}^2$ dans un voisinage
de $a$ dans $\mathbb{C}^2$. Supposons par la suite que $\cos
  t_i^*\not=\pm 1$; le cas o\`u par exemple $\cos t^*_1=1$, il suffit de
  remplacer $\cos (t_1-t_1^*)$ dans $\varphi$ par $\cos\sqrt{t_1-t_1^*}$.
L'application $\varphi$ est
un biholomorphisme local. Soient $g:=\varphi^{-1}\circ f \circ \varphi$
et $g_1:=\varphi^{-1}\circ f_1 \circ \varphi$. On a $g(t_1,t_2)=(dt_1,d
t_2)$ et $g\circ g_1= g_1\circ g$. Utilisant les d\'eveloppements de
Taylor de $g$ et $g_1$, on montre facilement que $g_1$ est
lin\'eaire. Posons $g_1(t_1,t_2)=(\alpha t_1+\beta t_2, \alpha'
t_1+\beta' t_2)$. Soit $f_1=(P,Q)$. On a $\cos(\alpha t_1'+\beta
t_2')=P(\cos t_1',\cos t_2')$ o\`u $t_i':=t_i-t_i^*$. 
Comme le membre \`a droite est p\'eriodique de
p\'eriode $(2\pi,2\pi)$, $\alpha$ et $\beta$ sont
entiers. L'\'egalit\'e pr\'ec\'edente implique $\sin\alpha
t_1'\sin\beta t_2'= \cos\alpha
t_1'\cos\beta t_2'-P(\cos t_1',\cos t_2')$. Si $\alpha\beta\not=0$,
$\sin\alpha t_1'$ s'\'ecrit en fonction de $\cos t_1'$. Ceci est
impossible car $\sin\alpha t_1'$ est impaire et non identiquement
nulle. Donc $\alpha\beta=0$. De m\^eme,
$\alpha'\beta'=0$. Comme $f_1$ et $g_1$ sont ouvertes, 
$\alpha=\beta'=0$ ou
$\alpha'=\beta=0$. 
On en d\'eduit facilement que
$P=\pm T_{d_1}z_1$ ou $\pm T_{d_1}z_2$. Le
polyn\^ome $Q$ l'est aussi. On conclut que $(f_1,f_2)$ se trouve dans
l'exemple 2.
\\
\underline{\it Cas 3.} 
Supposons que  $(f_1^m, f_2^n)$ se trouve dans l'exemple 3. Alors
  $f_1^m$ et $f_2^n$ sont homog\`enes. Par cons\'equent, $f_1$ et $f_2$
  sont homog\`enes et $(f_1,f_2)$ se trouve dans l'exemple 3.
\\
\underline{\it Cas 4.}
Supposons maintenant que $(f_1,f_2)$ se trouve dans l'exemple 
4. Soit $\pi$ l'application de $\mathbb{C}^2$ dans $\mathbb{C}^2$
avec $\pi(x,y)=(x+y,xy)$ ({\it voir} l'exemple 4). Quitte \`a remplacer
$(m,n)$ par $(2m,2n)$, on peut supposer que
$f_1^m\circ\pi=\pi\circ\tilde F_1$ et $f_2^n\circ\pi=\pi\circ\tilde
F_2$ o\`u $\tilde F_1(x,y):=(T_{d_1^m}x,T_{d_1^m}y)$,  $\tilde
F_2(x,y):=(T_{d_2^n}x,T_{d_2^n}y)$ ou $\tilde
  F_1(x,y)=(\lambda_1x^{d_1^m},\lambda_1 y^{d_1^m})$, $\tilde
  F_2(x,y)=(\lambda_2x^{d_2^n},\lambda_2 y^{d_2^n})$. 
Comme $f_1^m$ d\'efinit un rev\^etement d'un certain orbifold ${\cal
  O}=(\mathbb{P}^2,n)$ 
dans lui-m\^eme, $f_1$ d\'efinit \'egalement un rev\^etement de
${\cal O}$ 
dans lui-m\^eme. Remarquons que la courbe $H:=\{z_1^2-4z_2=0\}$ est
  invariante par $f_1^m$, $f_2^n$ et $n(H)=2$. 
Le fait que $f_i$ d\'efinit un rev\^etement de ${\cal O}$ dans
lui-m\^eme implique 
qu'il existe des
polyn\^omes $R_i$ v\'erifiant: 
$$f_{i1}^2-4f_{i2}=(z_1^2-4z_2)R_i^2.$$
Posons $H_i:=(H_{i1},H_{i2})=f_i\circ \pi$, on a: 
$$H_{i1}^2-4H_{i2}=(x-y)^2R_i^2(x+y,xy).$$
Ces \'egalit\'es signifient 
qu'il existe des applications $F_i$ v\'erifiant
$\pi\circ F_i=H_i$. 
On peut choisir les $F_i$ telles que $F_1\circ F_2=F_2\circ
F_1$, $F_1^m=\tilde F_1$ et $F_2^n=\tilde F_2$. 
Comme $F_1^m$, $F_2^n$ sont prolongeables
en des endomorphismes holomorphes de $\mathbb{P}^2$, $F_1$ et $F_2$ le sont
aussi. Comme dans les cas 2 et 3, 
on montre que 
$F_i=(\pm T_{d_i}x,\pm T_{d_i}y)$,  $(\pm T_{d_i}y,\pm
T_{d_i}x)$, $(\lambda_i x^{d_i},\lambda_i y^{d_i})$ ou $(\lambda_i
y^{d_i}, \lambda_i x^{d_i}$). Les deux premiers cas (resp. les deux
derniers) donnent les m\^emes applications
$f_i$. Donc $(f_1,f_2)$ appartient \`a l'exemple 4. 
\end{preuve}
\begin{proposition} Sous l'hypoth\`ese du th\'eor\`eme 1, 
pour toute composante $A$ de ${\cal
    C}_1\cup{\cal C}_2$, 
  il existe deux couples de nombres naturels 
  $(n, m)\not = (n',m')$ tels que
  $f_1^n f_2^m A=f_1^{n'}f_2^{m'} A$. De plus, il existe 
  $N$ qui ne d\'epend que du nombre des composantes de ${\cal
  C}_1\cup{\cal C}_2$
  tel que $n,m,n',m'$ soient major\'es par $N$.
\end{proposition}
\begin{remarque} Cette proposition reste valable pour les
    endomorphismes ouverts permutables d'une vari\'et\'e quelconque
    tels que les ensembles critiques des $f_i$ contiennent un nombre
    fini de composantes irr\'eductibles. Elle est d\'emontr\'ee en
    collaboration avec N. Sibony.
\end{remarque}
\begin{preuve} Il suffit de consid\'erer $A\subset{\cal C}_1$. 
La relation $f_1\circ f_2=f_2\circ
f_1$ implique:
$$f_1'(f_2(z)).f_2'(z)=f_2'(f_1(z)).f_1'(z).$$
D'o\`u:
\begin{equation}
\mathcal{C}_2\cup f_2^{-1}(\mathcal{C}_1)=
\mathcal{C}_1\cup f_1^{-1}(\mathcal{C}_2)
\end{equation}
Supposons qu'il n'existe pas $n,m,n',m'$ v\'erifiant cette proposition. 
Alors la suite $f_2^nA$ contient une infinit\'e
de courbes diff\'erentes. 
Ceci implique qu'il existe un $n$ tel que $f_2^nA\not
\subset \mathcal{C}_1$ car $\mathcal{C}_1$ ne contient qu'un nombre
fini de courbes irr\'eductibles. On choisit $n$ minimal
possible. Posons $A':=f_2^{n-1}A$.
La relation (3) implique que $A'\subset
\mathcal{C}_1\cap\mathcal{C}_2$ et $f_1^{-1}A'\subset\mathcal{C}_2\cup
f_2^{-1}\mathcal{C}_1$.
\par
Supposons que pour tout $m\geq 0$, $f_1^{-m}A'\subset \mathcal{C}_2$. Comme
$\mathcal{C}_2$ ne contient qu'un nombre fini de composantes, il existe
$m_1<m_2$ et une composante $A''\subset\mathcal{C}_2$ v\'erifiant 
$A'=f_1^{m_1}A''=f_1^{m_2}A''$. Par cons\'equent, $f_2^{n-1}A=
f_2^{n-1}f_1^{m_2-m_1}A$. C'est contradiction.
\par
Soit $m\geq 1$ le nombre minimal v\'erifiant $f_1^{-m}A'\not\subset
\mathcal{C}_2$. Comme $f_1^{-m+1}A'\subset {\cal C}_2$, 
(3) implique   $f_1^{-m}A'\subset\mathcal{C}_2\cup
f_2^{-1}\mathcal{C}_1$. Donc  $f_1^{-m}A'$ contient une composante de
$f_2^{-1}\mathcal{C}_1$. Il existe une composante $A_2$ de
$\mathcal{C}_1$ telle que $f_2A'=f_1^mA_2$. D'o\`u
$f_2^nA=f_1^mA_2$.
\par
De m\^eme, il existe une suite $\{A_i\}$ de composantes de
$\mathcal{C}_1$ avec $A_1=A$ 
et des $n_i$, $m_i\geq 1$ avec $n_1:=n$ et $m_1:=m$ tels que
$f_2^{n_i}A_i= f_1^{m_i}A_{i+1}$. Il existe  
$0<l<j$ tels que $A_l=A_j$. On d\'eduit que
$f_2^rf_1^s A=f_1^p A$ pour $r:=n_1+\cdots + n_{l-1}$,
$s:=m_l+\cdots + m_{j-1}$ et $p:=m_1+\cdots +m_{j-1}$. 
C'est contradiction.\par
Par le raisonnement ci-dessus, l'existence de $N$ est \'evidente. 
\end{preuve}
\section{Germes holomorphes et
ensembles critiques}
Soient $f$ une application ouverte d'un domaine $U\subset\mathbb{C}^n$
dans $\mathbb{C}^n$, $H$ et $H'$ deux hypersurfaces de $U$. 
Pour tout $a\in U$, on note $\mult(f,a)$ {\it la
multiplicit\'e} de $f$ en $a$, i.e. le nombre de pr\'eimages dans $U'$ 
d'un point $b$ g\'en\'erique suffisamment 
proche de $f(a)$, o\`u $U'$ est un voisinage
suffisamment petit de $a$. 
On note $\mult(f,H)$ {\it la multiplicit\'e} de $f$ sur $H$, i.e. la
multiplicit\'e de $f$ en un point g\'en\'erique de $H$ et
 $\mult(H\cap H',c)$ la multiplicit\'e de l'intersection
$H\cap H'$ en $c\in H\cap H'$. 
Posons $\m_f(H):= \mult(f,H)-1$. Alors $\m_f(H)>0$ si et seulement si 
$H$ appartient \`a l'ensemble critique de $f$ avec la multiplicit\'e
$\m_f(H)$.  Dans la suite, $f$ est un germe
d'application holomorphe, ouverte, d\'efinie au voisinage de
$0\in\mathbb{C}^2$ \`a l'image dans $\mathbb{C}^2$ et ${\cal C}$
son ensemble critique.
\begin{lemme} Supposons que
$f(x,y)=(x^d+\o(x^d),g(x,y))$ avec $d\geq 2$. Soient
${\cal C}_1$, ${\cal C}_2$, ..., ${\cal C}_m$
et $\{x=0\}$ les composantes irr\'eductibles de 
${\cal C}$. Posons $n_i:=\m_f({\cal C}_i)$ 
et $m_i:=\mult({\cal
  C}_i\cap\{x=0\}, 0)$. Alors 
$\sum m_in_i=\mult(h,0)-1=\m_h(0)$, o\`u $h:=f|_{\{x=0\}}$.
\end{lemme}
\begin{preuve} Soit $J$ le jacobien de $f$. Alors $\sum m_in_i$ est la
  multiplicit\'e en $0$ de la restriction de $|J|/x^{d-1}$ sur
  $\{x=0\}$. D'autre part,
  $$|J|= x^{d-1}\begin{array}{|cc|}
                d+\o(1) & \O(x) \\
                g_x    & g_y
                \end{array}
  $$
o\`u $g_x:=\partial g/\partial x$ et $g_y:=\partial g/\partial y$. Par
cons\'equent, 
$$\sum m_in_i=\mult(g_y(0,y),0)=\mult(g(0,y),0)-1=\mult(h,0)-1.$$
\end{preuve} 
\begin{lemme} Soient $f,h$ d\'efinis dans le lemme 3,
  $A$ une courbe lisse qui coupe $\{x=0\}$ transversalement en
  $0$. Soient $A_i$ les composantes irr\'eductibles de $f^{-1}A$,
  $n_i:=\mult(f,A_i)$ et $m_i:=\mult(A_i\cap\{x=0\},0)$. Alors
  $\mult(f,0)=d.\mult(h,0)=d.\sum m_in_i$.
\end{lemme}
\begin{preuve} Quitte \`a remplacer $f$ par $f\circ \sigma$, on peut
  supposer que $A=\{y=0\}$, o\`u $\sigma$ est une 
application holomorphe inversible.
Alors $0$ est une
solution de multiplicit\'e $\mult(f,0)$ de l'\'equation
$f(x,y)=0$. Posons $f=(l,g)$. La solution de l'\'equation $l=0$ est
$\{x=0\}$ avec la multiplicit\'e $d$. Comme $f^{-1}A=\{g=0\}$, on
a 
$$\mult(f,0)=d.\mult(\{g=0\}\cap\{x=0\},0)=d.\mult (h,0)=d\sum m_i
n_i.$$   
\end{preuve}
Dans la suite, on note $f_1$, $f_2$ deux germes d'application
holomorphe ouverte d\'efinies au voisinage de
$0\in\mathbb{C}^2$ \`a l'image dans $\mathbb{C}^2$ tels que
$f_1(0)=f_2(0)=0$ et $f_1\circ f_2=f_2\circ f_1$.
Notons ${\cal C}_i$ l'ensemble critique de $f_i$.
\begin{lemme} Supposons que 
$f_i(x,y)=(x^{d_1}+\o(x^{d_1}),a_iy+yh_i)$
o\`u $h_i(0)=0$, $a_i\not =0$ et $a_1^n\not= a_2^n$ pour tout
$n\in\mathbb{N}_+$. Soit ${\cal C}$ un germe de surface de Riemann
en $0$ (\'eventuellement singulier) tel que $f_1{\cal C}=f_2{\cal
  C}$. Alors
${\cal C}$ est l'une des deux courbes invariantes $\{x=0\}$ et $\{y=0\}$.
\end{lemme}
\begin{preuve}
Supposons que ${\cal C}\not =\{x=0\}$ et ${\cal C}\not =\{y=0\}$. 
On peut \'ecrire
${\cal C}=u(D)$
o\`u $D$ est un disque de $\mathbb{C}$ centr\'e en $0$,
$u(t)=(t^l,ct^m+\o(t^m))$ une
application holomorphe de $D$ dans $\mathbb{C}^2$ et $c\not = 0$. On a: 
$$f_i\circ u(t)=(t^{ld_i}+\o(t^{ld_i}),a_ict^m+\o(t^m)).$$ 
Comme $f_1{\cal C}=f_2{\cal C}$, on a $ld_1=ld_2$ et
$(a_1c)^{ld_1}=(a_2c)^{ld_2}$. D'o\`u $d_1=d_2$ et
$a_1^{ld_1}=a_2^{ld_1}$. C'est la contradiction recherch\'ee.
\end{preuve}
Soit $h(x,y)=\sum_{k,l\geq 1}a_{kl}x^ky^l$ 
une fonction holomorphe d\'efinie au voisinage de
$0\in\mathbb{C}^2$. Pour tout $\alpha\in\mathbb{R}_+$, on pose
$$d_\alpha:=\min_{a_{kl}\not=0}\alpha k+l.$$ 
Alors on peut \'ecrire $h=h_\alpha + \o_\alpha(y^{d_\alpha})$, o\`u 
$h_\alpha:=\sum_{\alpha k+l=d_\alpha} a_{kl} x^ky^l$   
et $\o_\alpha(y^{d_\alpha})$ est une fonction v\'erifiant
$\o_\alpha(y^{d_\alpha})=\o(t^{d_\alpha})$ quand $x=\O(t^\alpha)$ et
$y=\O(t)$.
Consid\'erons $h=ay^d+xg(x,y)$ avec $a\not
=0$ et $g$ holomorphe. On d\'efinit un nombre rationnel $\alpha_h\geq 1$ par: 
$$\alpha_h:=\max\{1, \min_{d_\alpha=d}\alpha\}.$$
\begin{proposition}  Soient
$f_i=(x^{d_i}+\o(x^{d_i}),h_i)$ et
$\alpha:=\max(\alpha_{h_1},\alpha_{h_2})$
o\`u
$h_i=y^{d_i}+xg_i(x,y)$. Supposons que
$d_1^m\not=d_2^n$ pour tout $(m,n)\in\mathbb{N}^2-
(0,0)$. 
Soient 
$P_i(y):=h_{i,\alpha}(1,y)$. Alors $P_1\circ
P_2=P_2\circ P_1$ et l'une
des conditions suivantes est vraie:
\begin{enumerate}
\item $\alpha=1$, $(P_1,P_2)$ est conjugu\'e \`a $(y^{d_1},\gamma y^{d_2})$
  avec $\gamma^{d_1-1}=1$;
\item $\alpha=1$, $(P_1,P_2)$ est conjugu\'e \`a $(\pm \tilde
  T_{d_1},\pm \tilde T_{d_2})$; 
\item $\alpha=2$, $\beta^{-1}P_i(\beta y) =\pm
  \tilde T_{d_i}(y)$ o\`u $\beta^{d_i-1}=1$.
\end{enumerate}
\end{proposition}
\begin{preuve} 
Si $\alpha=1$, $[x^{d_i}:h_{i,\alpha}]$  
d\'efinissent deux
endomorphismes polynomiaux permutables de $\mathbb{P}^1$ dont $[x:y]$
sont les coordonn\'ees homog\`enes. D'o\`u $P_1\circ P_2=P_2\circ P_1$.
Comme $d_1^m\not=d_2^n$ pour $(m,n)\not= (0,0)$, $(P,Q)$
est conjugu\'e \`a $(y^{d_1},\gamma y^{d_2})$ ou \`a $(\pm \tilde T_{d_1},
\pm \tilde T_{d_2})$. 
L'une des conditions 1 ou 2 est vraie.
\par 
Supposons maintenant que $\alpha=p/q$ avec  $p$, $q$ premiers
entre eux et $p>q\geq 1$. 
La relation $f_1\circ f_2=f_2\circ f_1$ appliqu\'ee \`a
$x=t^p$ et $y=at^q$ nous donne $P_1\circ P_2= P_2\circ
 P_1$. Comme
$d_1^m\not=d_2^n$ pour $(m,n)\not =(0,0)$, $(P_1,P_2)$ est
conjugu\'e \`a $(a^{d_1},\gamma a^{d_2})$ ou \`a $(\pm \tilde T_{d_1},\pm
\tilde T_{d_2})$. 
\par
Le coefficient de $a^{d_i-1}$ dans $P_i(a)$ est nul car $\alpha>1$. 
Alors si $(P_1,P_2)$ est
conjugu\'e \`a $(a^{d_1},\gamma a^{d_2})$, il sera \'egal \`a
$(a^{d_1}, \gamma' a^{d_2})$. 
Ceci est impossible ({\it voir} la d\'efinition de
$\alpha=\max(\alpha_{h_1},\alpha_{h_2})$). 
Alors $(P_1,P_2)$ est conjugu\'e \`a $(\pm
\tilde T_{d_1},\pm \tilde T_{d_2})$.
\par
Par la d\'efinition de $\alpha$, $ P_i(a)=a^mP_i^*(a^p)$ pour un certain
polyn\^ome $P^*_i$ . Alors l'ensemble des z\'eros de $P_i$ est
stable par la rotation $a\longmapsto \sqrt[p]{1}a$. On en d\'eduit que
$p=2$, donc $q=1$ et $\alpha=2$ car les z\'eros de $\pm
\tilde T_{d_i}$ sont align\'es.
\par
Finalement, il existe une application lin\'eaire $\sigma(a):=\beta
a+\theta$ avec $\beta\not =0$ telle que $\sigma^{-1}\circ
P_i\circ\sigma=\pm \tilde T_{d_i}$. 
Si $d_i$ est pair (ou impair) les fonctions
$\tilde T_{d_i}$ et $P_i$ le sont aussi. De plus, le coefficient de $y^{d_i-1}$
dans $P_i(y)$ est nul; ceux de $y^{d_i}$ dans $P_i$ et dans
$\tilde T_{d_i}$ sont 1. Donc $\theta=0$ et $\beta^{d_i-1}=1$.
\end{preuve}
\begin{corollaire} Si la condition 3 de la proposition 2 est vraie,
alors il existe une  
courbe ${\cal D}=\{x=\beta^2y^2/4+\o(y^2)\}$ invariante par $f_1$ et $f_2$ 
  telle que $\overline{f_i^{-1}{\cal D}\setminus{\cal D}}={\cal
    C}_i$. De plus,
\item[\ \ ] 
{\rm 1.} Si $d_i$ est impair, $\mathcal{C}_i$ est une r\'eunion de $(d_i-1)/2$
  courbes ${\cal D}_r=\{x=\beta^2y^2/r+\o(y^2)\}$ o\`u
  $r\in\mathbb{R}^+$ et $\pm\sqrt{r}\in
  \tilde T_{d_i}^{-1}\{\pm 2\}\setminus \{\pm 2\}$. 
\item[\ \ ] 
{\rm 2.} Si $d_i$ est pair, $\mathcal{C}_i$ est une r\'eunion de $d_i/2-1$
  courbes ${\cal D}_r=\{x=\beta^2y^2/r+\o(y^2)\}$ et d'une courbe
  ${\cal D}_0$ qui coupe $\{x=0\}$ transversalement o\`u
  $r\in\mathbb{R}^+$ et $\pm\sqrt{r}\in
  \tilde T_{d_i}^{-1} \{\pm 2\}\setminus \{\pm 2,0\}$. 
\end{corollaire}
\begin{preuve} On peut supposer que $i=1$. Quitte au changement de
  coordonn\'ees $(z_1,z_2)\mapsto (z_1,\beta z_2)$, on peut supposer
  que $\beta=1$.   
On \'ecrit 
$$f_j=(x^{d_j}+\o(x^{d_j}),
 \pm x^{d_j/2}\tilde T_{d_j}(y/\sqrt{x})+\o_2(y^{d_j})).$$ 
Soit $\tau$ l'application
  de $\mathbb{C}^2$ dans $\mathbb{C}^2$ d\'efinie
  par $\tau(t,a):=(t^2,at)$. Cette application est localement
  inversible sauf sur $\{t=0\}$. On consid\`ere
  $F_j:=\tau^{-1}\circ f_j\circ\tau$. Avec les $f_j$ d\'ecrits
ci-dessus, on peut trouver 
  des $F_j$ holomorphes au voisinage de $\{t=0\}$ v\'erifiant 
$F_1\circ F_2=F_2\circ F_1$ et 
$$F_j(t,a)=(t^{d_j}+\o(t^{d_j}),\pm \tilde T_{d_j}(a)+\O(t)).$$
Soit ${\cal C}'_j\cup\{t=0\}$
l'ensemble critique de $F_j$. On a $\tau{\cal C}'_j={\cal
  C}_j$. Fixons un $N$ et un $R$ suffisament grands. Il existe un
$\epsilon>0$ tel que les applications $F_j^l$ soient d\'efinies sur
$\{|t|<\epsilon\}\times \{|a|<R\}$ pour tout $l=1,2,\ldots, N$.
On consid\`ere ici le cas o\`u $d_1$ est pair et les signes $\pm$ sont
les $+$; les autres cas  seront
trait\'es de m\^eme mani\`ere.
\par
Consid\'erons un $b\in \tilde T_{d_1}^{-1}\{2\}\setminus\{0,2\}$. 
Alors $b$ est un
point critique d'ordre 1
de $\tilde T_{d_1}$. D'apr\`es le lemme 3, il existe
une courbe ${\cal D}'_b\subset{\cal C}_j'$ qui coupe $\{t=0\}$
transversalement en $(0,b)$. 
La proposition 1 s'applique \'egalement pour les $F_j$. 
Il existe $(m,n)\not=(m',n')$ tels que $F_1^mF_2^n{\cal
  D}'_b=F_1^{m'}F_2^{n'} {\cal D}'_b$ et $m,n,m',n'$ soient plus petits
que $N-1$. On en d\'eduit que $F_1^mF_2^n(F_1{\cal
  D}'_b)=F_1^{m'}F_2^{n'} (F_1{\cal D}'_b)$. Or $F_1{\cal D}'_b$ passe
par $(0,2)$ un point fixe de $F_1$ et $F_2$, 
le lemme 5 (appliqu\'e \`a $F_1^mF_2^n$ et $F_1^{m'}F_2^{n'}$)  
montre que $F_1{\cal D}_b$ est invariante
par $F_1$ et $F_2$. Cette courbe invariante, not\'ee ${\cal D}'$, est
la courbe stable passant par $(0,2)$; elle  
est unique, lisse et ind\'ependante de $b$; 
elle coupe $\{t=0\}$ transversalement. La courbe ${\cal D}'$ 
est d\'efinie par une
\'equation du type $a=2+\O(t)$. Ceci implique que ${\cal
  D}:=\tau^{-1}{\cal D}'$ est d\'efinie par $y^2=4x+x\tilde\o(1)$ o\`u
la notation
$\tilde\o(1)$ 
signifie une fonction multi-valante tendant vers $0$ quand
$x\rightarrow 0$. On a $\mult({\cal
  D}\cap\{x=0\})=2$ et $\{x=0\}$ est tangente \`a ${\cal D}$. 
Donc ${\cal D}$ est lisse, irr\'eductible, tangente \`a
$\{x=0\}$. Elle est invariante par $f_1$, $f_2$ et 
d\'efinie par une \'equation du type $x=y^2/4+\o(y^2)$.
Posons $r:=b^2$, ${\cal D}_r:=\tau^{-1}{\cal D}'_b$. De m\^eme
mani\`ere, on prouve que ${\cal
  D}_r=\{x=y^2/r+\o(y^2)\}$. On a ${\cal D}_r\subset{\cal C}_1$. 
Il y a $d_1/2-1$ tels nombres $r$ diff\'erents. 
Comme $\mult
(f_1|_{\{x=0\}},0)=d_1$, le lemme 3 implique que ${\cal C}_1$ est la
r\'eunion des ${\cal D}_r$ et d'une courbe ${\cal D}_0$ qui coupe
$\{x=0\}$ transversalement en $0$. Remarquons qu'il existe 
une composante ${\cal
  D}''$ de $\tau^{-1} (f_1{\cal  D}_0)$ passant par $(0,-2)$ et 
$(n,m)\not=(n',m')$ tels que $F_1^nF_2^m {\cal
  D}''=F_1^{n'}F_2^{m'}{\cal D}''$. D'o\`u $F_1{\cal D}''={\cal
  D}'$. Ceci montre aussi que ${\cal D}''$ est la composante de
$\tau^{-1}({\cal D})$ qui passe par $(0,-2)$. D'o\`u
${\cal D}_0\subset f_1^{-1}{\cal D}$.  
 Le lemme 4 montre que
$f_1^{-1}{\cal D}$ contient seulement les ${\cal D}_r$, ${\cal D}_0$ et
${\cal D}$.
\end{preuve}
\section{Preuve du th\'eor\`eme 1: premier cas} 
Dans ce paragraphe, on suppose que 
$f_i|_\infty=[w_1^{d_i}:w_2^{d_i}]$ ({\it voir} l'introduction).
D'apr\`es le lemme 1, on a $(f_1|_\infty)^n\not= (f_2|_\infty)^m$
pour tout $(n,m)\not=(0,0)$. Donc $d_1^n\not=d_2^m$ pour tout
$(n,m)\not =(0,0)$.
Soient $a=[0:0:1]$ et
  $a'=[0:1:0]$. En utilisant les coordonn\'ees locales
  $x:=w_0/w_2$, $y:=w_1/w_2$, au voisinage de $a$, on peut \'ecrire
$f_i(x,y)=(x^{d_i}+\o(x^{d_i}), h_i)$
o\`u $h_i=y^{d_i}+xg_i$.
\par
D'apr\`es la
proposition 2, il existe une application lin\'eaire $\sigma(z_1)=\beta
z_1+\theta$ telle que l'une des conditions suivantes soit vraie:
\begin{enumerate}
\item   
$\sigma^{-1}\circ f_{11}\circ\sigma=z_1^{d_1}$ et 
$\sigma^{-1}\circ f_{21}\circ\sigma=\gamma z_1^{d_2}$ avec
$\gamma^{d_1-1} =1$;
\item  
$\sigma^{-1}\circ f_{i1}\circ\sigma=\pm \tilde T_{d_i}z_1$;
\item $f_{i1}=\pm
  \beta^{-1} 
z_2^{d_i/2}\tilde T_{d_i}(\beta z_1/\sqrt{z_2})+\Xi_i$ o\`u 
$\Xi_i=\sum_{k+2l<d_i}b_{ikl}z_1^kz_2^l$.
\end{enumerate} 
Quitte \`a remplacer $f_i$ par $f_i^2$, on peut supposer que les signes
$\pm$ sont les $+$. Comme $h_i=y^{d_i}+xg_i$, dans les
conditions 1, 2, on a $\beta^{d_i-1}=1$. 
Quitte \`a effectuer un changement de coordonn\'ees
$(z_1,z_2)\mapsto (z_1+\theta/\beta,z_2)$ dans les conditions 1, 2, 
on peut supposer que $\theta=0$. Alors dans 1., on peut supposer que
$\sigma=\Id$; dans 2. et 3.,  
en utilisant le changement $(z_1,z_2)\mapsto
(\beta z_1,z_2)$, on peut \'egalement
supposer que $\beta=1$ et $\sigma=\Id$.
De m\^eme pour $a'$, l'une des conditions suivantes est vraie:
\begin{enumerate}
\item[1'.] $f_{i2}=z_2^{d_i}$;
\item[2'.] $f_{i2}=\tilde T_{d_i}z_2$;
\item[3'.] $f_{i2}=
   {\beta'}^{-1}z_1^{d_i/2}
\tilde T_{d_i}(\beta' z_2/\sqrt{z_1})+\Xi'_i$ o\`u 
$\Xi_i'=\sum_{2k+l<d_i}b'_{ikl}z_1^kz_2^l$, ${\beta'}^{d_i-1}=1$.
\end{enumerate} 
Il est clair que si les conditions 3 et 3' sont fausses, le couple
$(f_1,f_2)$ se trouve dans les exemples 1-3 \`a une classe de
conjugaison pr\`es. La preuve du th\'eor\`eme 1 sera complet\'ee par
les propositions qui suivent.
\begin{proposition} Si les conditions 1, 3' ou 1',3 sont
  r\'ealis\'ees, le couple $(f_1,f_2)$ se trouve dans l'exemple 4 \`a
  une classe de conjugaison pr\`es. 
\end{proposition}
\begin{preuve} On peut supposer que 1' et 3 sont vraies. 
D'apr\`es la proposition 2, il existe une courbe
  alg\'ebrique ${\cal D}$  invariante tangente \`a $\{w_0=0\}$ en
  $a$. Comme $f_i|_\infty=[w_1^{d_i}:w_2^{d_i}]$ et
  comme $f_{i2}=z_2^{d_i}$, on a
${\cal D}\cap\{w_0=0\}=\{a\}$; ${\cal D}$ est donc
  de degr\'e 2. L'ensemble exceptionnel de $f_i|_{\cal D}$ qui est de
  cardinal $\leq 2$ et qui contient les points
  d'intersection de ${\cal D}$ avec $\{w_0=0\}\cup\{w_2=0\}$. Par
  cons\'equent, ${\cal D}$ est tangente \`a $\{w_2=0\}$ en un point, not\'e
  $b$. Quitte \`a effectuer un changement de coordonn\'ees du type 
$z_2\mapsto z_2+\alpha$, on peut supposer que $b=0$. 
Alors l'application $f_i|_{\{w_1=0\}}$  poss\`ede deux points
  exceptionnels; elle est donc \'egale \`a
  $z_2^{d_i}$. Appliquant la proposition 2 en $b$,   
on constate que
  $\Xi_i=0$. On v\'erifie facilement que les 
  $f_i=(z_1^{d_i/2}\tilde 
  T_{d_i}(z_1/\sqrt{z_2}), z_2^{d_i})$ sont conjugu\'ees aux
  applications construites dans l'exemple 4 pour $h_i(x)=x^{d_i}$.
\end{preuve}
\begin{proposition} Les conditions 2, 3' (ou 2', 3) ne sont pas
  r\'ealis\'ees simultanement.
\end{proposition}
\begin{preuve} Supposons, par exemple, que 2 et 3' sont
  vraies. 
Soit ${\cal D}$ la courbe invariante de $f_i$ passant par $a'$ qui est
  d\'ecrite dans le corollaire 2. Cette courbe est tangente \`a 
$\{w_0=0\}$ en un
  point unique $a'$. Par sa description dans le corollaire 2, ${\cal D}$
  est de degr\'e 2. L'application $f_1$ agit sur ${\cal D}$
comme un polyn\^ome de degr\'e $d_1$. Alors il y a au plus
$d_1-1$ points critiques de $f_1|_{\cal D}$ qui sont diff\'erents
de $a'$.
\par 
Posons ${\cal E}_s=\{z_1=s\}$.
Observons que 
  $I:=\tilde T_{d_1}^{-1}\{\pm 2\}\setminus \{\pm 2\}$ est
l'ensemble critique de $\tilde T_{d_1}$. Comme $f_{11}=\tilde
T_{d_1}(z_1)$, si $s\in I$ et
si ${\cal E}_s$ coupe ${\cal D}$ transversalement
en $z$ alors $z$ est un point critique de $f_1|_{\cal D}$. De
plus, il y a une seule droite ${\cal E}_s$ qui ne coupe pas
${\cal D}$ transversalement en deux points car ${\cal D}$ est une courbe
rationnelle ne passant pas par $[0:0:1]$. Alors $f_1|_{\cal D}$
poss\`ede au moins $2(d_1-2)$ points critiques diff\'erents de
$a'$ car $\#I=d_1-1$. C'est la contradiction recherch\'ee.
\end{preuve}
\begin{proposition} Les conditions 3, 3' ne sont pas
  r\'ealis\'ees simultanement.
\end{proposition}
\begin{preuve} Supposons que 3, 3' sont vraies. Pour simplifier les
  calculs on suppose que $\beta'=1$. Quitte \`a remplacer $f_i$
par l'un de ses it\'et\'es, on peut supposer que  $d_i>>1$.
D'apr\`es le
  corrollaire 2, il existe une courbe invariante ${\cal D}$ tangente
  \`a $\{w_0=0\}$ en exactement deux points $a$ et
  $a'$ telle que $f_i^{-1}{\cal D}={\cal D}\cup
{\cal C}_i$. 
Cette courbe ${\cal D}$ est de degr\'e 4.
D'apr\`es les conditions 3 et 3', quitte \`a effectuer un changement
  de coordonn\'ees du type $(z_1,z_2)\mapsto (z_1+\alpha,z_2+\beta)$,
  on peut supposer que ${\cal D}$ est d\'efinie par $\Phi=0$ o\`u 
$$\Phi=z_1^2z_2^2-4z_1^3-4z_2^3+rz_1^2+uz_1z_2+vz_2^2+\o(|z|^2)$$
est un polyn\^ome de degr\'e 4. D'apr\`es le corollaire 2, il existe
un polyn\^ome $W$ tel que $\Phi(f_1)=\Phi.W^2$. 
\par 
Soit $\Delta_i$ la
somme des termes $z_1^mz_2^n$ de $f_{11}$ qui v\'erifient
$m+2n=d_1-i$. Posons $\Theta_i(\alpha):=\Delta_i(\alpha,1)$.
Montrons que $\Theta_1=\Theta_2=0$. Posons $z_1= \alpha t$ et
$z_2=t^2$. Notons $W_i$ la somme des termes $z_1^mz_2^n$ de $W$
qui v\'erifient $m+2n=3d_1-3-i$. Posons
$V_i(\alpha):=W_i(\alpha,1)$. 
Alors d'apr\`es les conditions 3, 3', on a 
$f_{11}(\alpha t,t^2)=\tilde
T_{d_1}(\alpha)t^{d_1}+\Theta_1(\alpha)t^{d_1-1}+\o(t^{d_1-1})$,
$f_{12}(\alpha t,t^2)= t^{2d_1}+\o(t^{2d_1-1})$. On peut \'ecrire
$W(at,t^2)=V_0(\alpha)t^{3d_1-3}+V_1(\alpha)t^{3d_1-4}+ \o(t^{3d_1-4})$.
En identifiant les coefficients de $t^{6d_1}$ et $t^{6d_1-1}$
dans l'\'equation  $\Phi(f_1)=\Phi W^2$, 
on obtient $\tilde T_{d_1}^2(\alpha)-4=(\alpha^2-4)V_0^2(\alpha)$ et 
$2\
\tilde T_{d_1}(\alpha)\Theta_1(\alpha,1)
=2(\alpha^2-4)V_0(\alpha)V_1(\alpha)$. 
La relation $\tilde
T_{d_1}(\alpha)^2-4=(\alpha^2-4)V_0^2(\alpha)$ implique que $V_0$
est la d\'eriv\'ee de $\tilde T_{d_1}$ et que 
$\tilde T_{d_1}(\alpha)$,
$(\alpha^2-4)V_0(\alpha)$ sont premiers entre eux. Alors $V_0$
est de degr\'e $d_1-1$ et 
le polyn\^ome
$\Theta_1(\alpha)$ est divisible par $2(\alpha^2-4)V_0(\alpha)$.
On en
d\'eduit que $\Theta_1=0$ car $\deg
\Theta_1<d_1$. On a \'egalement $V_1=0$. Alors $f_{11}$ ne
contient aucun terme $z_1^mz_2^n$ avec $m+2n=d_1-1$. De m\^eme
mani\`ere, on montre que  $f_{12}$ ne
contient aucun terme $z_1^mz_2^n$ avec $2m+n=d_1-1$. Alors on peut
\'ecrire $f_{11}(\alpha t,t^2)=\tilde
T_{d_1}(\alpha)t^{d_1}+\Theta_2(\alpha)t^{d_1-2}+\o(t^{d_1-2})$, 
$f_{12}(\alpha t,t^2)=t^{2d_1}+\o(t^{2d_1-2})$ et
$W=V_0(\alpha)t^{3d_1-3}+V_2(\alpha)t^{3d_1-5}+ \o(t^{3d_1-5})$.
En identifiant les coefficients de $t^{6d_1-2}$
dans l'\'equation  $\Phi(f_1)=\Phi
W^2$, on obtient
$2\tilde
T_{d_1}(\alpha)\Theta_2(\alpha)=2V_0(\alpha)[(\alpha^2-4)V_2(\alpha)+
vV_0(\alpha)]$. 
Alors $\Theta_2(\alpha)$ est divisible par $V_0(\alpha)$ et donc
$\Theta_2=0$. 
\par
On peut donc \'ecrire $f_{11}(\alpha t,t^2)=\tilde
T_{d_1}(\alpha)t^{d_1}+\Theta_3(\alpha)t^{d_1-3}+\o(t^{d_1-3})$.
Comme le coefficient de $z_1z_2^{d_1-2}$ dans $f_{12}$ est \'egal \`a
$-d_1$, on a 
$f_{12}(\alpha t,t^2)=t^{2d_1}-d_1\alpha
t^{2d_1-3}+\o(t^{2d_1-3})$. On a \'egalement
$W=V_0(\alpha)t^{3d_1-3}+V_2(\alpha)t^{3d_1-5} +
V_3(\alpha)t^{3d_1-6}+ \o(t^{3d_1-6})$.
En identifiant les coefficients de $t^{6d_1-3}$ dans l'\'equation
 $\Phi(f_1)=\Phi W^2$, on obtient
$2\tilde T_{d_1}(\alpha)
[\Theta_3(\alpha)-d_1\alpha \tilde T_{d_1}(\alpha)]=
2V_0(\alpha)[(\alpha^2-4)V_3(\alpha)+ uV_0(\alpha)]$.
Par cons\'equent, 
$\Theta_3(\alpha)-d_1\alpha \tilde
T_{d_1}(\alpha)$ est divisible par $V_0(\alpha)=\tilde
T_{d_1}'(\alpha)$.
\par
L'application $f_1^2$ doit v\'erifier une propri\'et\'e analogue.
Par les calculs simples, on constate que $\tilde T_{d_1}'\circ\tilde
T_{d_1}(\alpha) \Theta_3(\alpha)-\frac{3}{2}d_1^2\alpha \tilde
T_{d_1^2}(\alpha) +\frac{1}{2}d_1a^2\tilde T_{d_1}'\circ\tilde T_{d_1}(\alpha)$ est
divisible par $\tilde T_{d_1^2}'(\alpha)$.
Comme $\tilde T_{d_1^2}'$ est divisible par $\tilde T_{d_1}'
\circ\tilde T_{d_1}$, le polyn\^ome
$\alpha\tilde T_{d_1^2}(\alpha)$ est divisible par
$\tilde T_{d_1}'\circ\tilde T_{d_1}(\alpha)$.
C'est une contradiction car
$\tilde T_{d_1^2}$ et sa d\'eriv\'ee sont premiers entre eux.
\end{preuve}
\section{Preuve du th\'eor\`eme 1: deuxi\`eme cas} 
Dans ce paragraphe, on suppose que 
$f_i|_\infty=[w_1^{d_i}:w_1^{d_i} \tilde 
T_{d_i}(w_2/w_1)]$ ({\it voir} l'introduction). 
D'apr\`es le lemme 1, on a $(f_1|_\infty)^n\not= (f_2|_\infty)^m$
pour tout $(n,m)\not=(0,0)$. Donc $d_1^n\not=d_2^m$ pour tout
$(n,m)\not =(0,0)$.
Comme dans le paragraphe
  pr\'ec\'edent, on utilise la proposition 2 et le corollaire 2
  en $a:=[0:0:1]$. Quitte \`a remplacer $f_i$ par $f_i^n$  et \`a
  changer les  coordonn\'ees,  
l'une des conditions suivantes sera vraie:
\begin{enumerate}
\item 
$f_{i1}=z_1^{d_i}$;
\item 
$f_{i1}=\tilde T_{d_i}z_1$;
\item $f_{i1}=
z_2^{d_i/2}\tilde T_{d_i}(z_1/\sqrt{z_2})+\Xi_i$ o\`u 
$\Xi_i=\sum_{k+2l<d_i}b_{ikl}z_1^kz_2^l$.
\end{enumerate} 
Soient $b:=[0:1:2]$ et $c:=[0:1:-2]$.
\begin{lemme} Supposons que 
  $d_i$ est impair. Il existe une courbe ${\cal D}$ 
  invariante par $f_1$ et $f_2$, passant par $a$, $b$ et
  $c$ telle que $\mult({\cal D}\cap\{w_0=0\},a)=2$, $\mult({\cal
  D}\cap \{w_0=0\},b)=\mult({\cal D}\cap \{w_0=0\},c)=1$
et $f_i^{-1}{\cal D}={\cal D}\cup{\cal C}_i$.
De plus, on a $\mult(f_i,{\cal C}_i)=2$ et $\mult({\cal
  C}_i\cap\{w_0=0\},d)=1$ pour tout $d\in {\cal C}_i\cap\{w_0=0\}$ et
  $d\not =a$.
\end{lemme}
\begin{preuve} On peut supposer que $i=1$.  
Soient $U$ un voisinage suffisamment petit de $\{w_0=0\}$
  et $A_p$ les composantes irr\'eductibles de ${\cal C}_1\cap U$ qui
  ne passent pas par $a$. Remarquons que $b$ et $c$ sont fixes et que
  l'ensemble critique de $f_1|_\infty$ est \'egal \`a $I\cup\{a\}$ o\`u
  $I:=(f_1|_\infty)^{-1} \{b,c\}\setminus \{b,c\}$.
D'apr\`es le lemme 3, $A_p$ coupe $\{w_0=0\}$ transversalement 
en un point de $I$ et
  $\mult(f_1,A_p)=2$. D'apr\`es la proposition 1, 
  il existe $(n,m)\not=(n',m')$ tels que
  $f_1^nf_2^mA_p=f_1^{n'}f_2^{m'}A_p$. Ceci implique
  $f_1^nf_2^m(f_1A_p)=f_1^{n'}f_2^{m'}(f_1A_p)$. D'apr\`es le lemme 5,
  $f_1A_p$ est la vari\'et\'e stable de $f_1$ et $f_2$ en $b$ ou $c$. Posons
  ${\cal D}'$ la courbe alg\'ebrique 
  invariante par $f_1$ et $f_2$ qui contient les
  $f_1A_p$.
\par
Si la condition 1 ci-dessus est vraie, ${\cal D}'$ ne passe pas par
  $a$. Il suffit de poser ${\cal D}:={\cal D}'\cup\{w_1=0\}$.
\par
Si la condition 2 ci-dessus est vraie, ${\cal D}'$ ne passe pas par
  $a$. Il suffit de poser 
${\cal D}:={\cal D}'\cup \{w_1=2\}\cup\{w_1=-2\}$.
\par
Lorsque la condition 3 ci-dessus est vraie,
on pose ${\cal D}:={\cal D}'$ si
  ${\cal D}'$ passe par $a$, sinon ${\cal D}$ est la r\'eunion de
  ${\cal D}'$ avec la courbe alg\'ebrique contenant la courbe
  invariante passant par $a$ qui est d\'ecrite dans
le corollaire 2. 
On a $\overline{f_1^{-1}{\cal D}\setminus{\cal D}}\subset{\cal C}_1$
  et $\deg \overline{f_1^{-1}{\cal D}\setminus{\cal
  D}}=2d_1-2\geq\deg{\cal C}_1$. D'o\`u
  $\overline{f_1^{-1}{\cal D}\setminus{\cal D}}= {\cal C}_1$ et
  $\mult(f_1,{\cal C}_1)=2$. 
\end{preuve}
\begin{lemme} Supposons que 
$d_i$ est pair. Il existe une droite $L$ 
  invariante par $f_1$, $f_2$, passant par $b$. La courbe $\overline{
f_i^{-1}L\setminus L}$ contient une droite $L'$ passant
  par $c$. 
La r\'eunion des composantes de ${\cal C}_i$
qui ne passent pas par $a$,
  est \'egale \`a $H:=\overline{f_i^{-1}(L\cup L')
\setminus (L\cup L')}$.
De plus, $H$ coupe $\{w_0=0\}$ transversalement  et $\mult(f_i,H)=2$.
\end{lemme}
\begin{preuve} On peut supposer que $i=1$. 
On remarque que $b=f(c)$ est un point  fixe et que
  l'ensemble critique de $f_1|_\infty$ est \'egal \`a $I\cup\{a\}$
  o\`u $I:=(f_1|_\infty)^{-1} \{b,c\}\setminus \{b,c\}$.
On utilise les  notations
du lemme pr\'ec\'edent. Soit $A_p$ une composante telle que
  $f_1A_p$ passe par $b$. On d\'emontre comme dans la proposition
  pr\'ec\'edente que $A:=f_1A_p$ est invariante. D'apr\`es les lemmes
  3 et 4,  $f_1^{-1}A$ contient une courbe $B\not\subset {\cal C}_1$ 
  qui coupe $\{w_0=0\}$
  transversalement en $c$. Ceci implique que la courbe alg\'ebrique
  irr\'eductible, contenant $A$ (not\'ee $L$) 
ne peut pas passer par $a$ car sinon
  d'apr\`es le corollaire 2, $\overline{f_1^{-1}L\setminus L}\subset
  {\cal C}_1$. Alors $L$ coupe $\{w_0=0\}$
  transversalement en un seul point.
Elle est donc une droite. La courbe
  alg\'ebrique irr\'eductible $L'$ qui contient $B$ coupe
  $\{w_0=0\}$ transversalement en un point unique $c$. Elle est
  \'egalement une droite. 
  Utiliser les lemmes 3 et 4,
on montre que $H\subset {\cal C}_1$, $H$
  coupe $\{w_0=0\}$ transversalement et
  $\mult(f_i,H)=2$.  
\end{preuve}
\begin{proposition}
Si la condition 3 est vraie, le couple $(f_1,f_2)$ se
  trouve dans l'exemple 4 \`a une classe de conjugaison pr\`es.
\end{proposition}
\begin{preuve}
On consid\`ere le cas o\`u $d_1$ est impair. Le cas contraire sera
  trait\'e de m\^eme mani\`ere. Soit ${\cal D}$ la courbe d\'efinie
  dans le lemme 6. C'est une courbe invariante de degr\'e 4. Comme
  $f_1|_{\cal D}$ poss\`ede au plus deux points exceptionnels, ${\cal
  D}$ est r\'eductible. Il y a donc deux cas possibles.
\par
Supposons que ${\cal D}$ est la r\'eunion d'une courbe irr\'eductible
${\cal E}$ de degr\'e 3 et d'une 
droite $L$. Alors $f_1|_{\cal E}$ est conjugu\'e \`a $x^{\pm d_1}$ car
il poss\`ede deux points exceptionnels (qui se trouvent dans $\{w_0=0\}$). 
Notons $I:={\cal E}\cap
L$. Alors $f_1^{-1}L$ passe par $J:=(f_1|_{\cal E})^{-1}I$. Comme $\deg
f_1^{-1}L=(d_1+1)/2$, on a $\# J\leq 3(d_1+1)/2$. On en d\'eduit que
$\#I=1$ car $f_1|_{\cal E}$ est conjugu\'e \`a $x^{\pm d_1}$. Soit $d$
le point d'intersection de ${\cal E}$ et $L$. C'est un point fixe de
$f_1$. La valeur propre de $f_1$ en $d$ suivante la direction tangente \`a
${\cal E}$ est \'egale \`a $d_1$. On montre facilement que la valeur
propre de $f_1$ en $d$ suivante la direction $L$ est diff\'erente de $0$ et
de $\pm d_1^2$.
Les applications $f_i$ agissent sur
$L$ comme les polyn\^omes permutables. Alors $f_1|_L$ est conjugu\'e
\`a $x^{d_1}$ ou \`a $T_{d_1}$. Le fait que la valeur propre de
$f_1|_L$ en $d$ est diff\'erente de $0$ et de $\pm d_1^2$ implique que
$f_1|_L$ est r\'egulier en chaque point de $K:=(f_1|_L)^{-1}d\setminus
\{d\}$ et que $\# K=d_1-1$. L'ensemble $f_1^{-1}{\cal E}\setminus
{\cal E}$
passe par $K$. Alors $K$ appartient \`a ${\cal C}_1$. On en d\'eduit
que la courbe $\overline{f_1^{-1}L\setminus L}$ passe par $K$. C'est une
contradiction car le degr\'e de cette courbe est \'egal \`a $(d_1-1)/2<\#K$. 
\par
Alors la courbe
${\cal D}$ contient une courbe rationnelle $R$ tangente \`a
$\{w_0=0\}$ en $a$. D'apr\`es le corollaire 2, quitte \`a changer
des coordonn\'ees, on peut supposer que  $R=\{z_1^2-4z_2=0\}$. 
D'apr\`es le corollaire 2, appliqu\'e au point $a$, on a 
$\overline{f_i^{-1}R\setminus R}\subset{\cal C}_i$. 
De plus, la multiplicit\'e de $f_i$ sur cet ensemble est
\'egale \`a 2. Par cons\'equent, il existe des polyn\^omes $R_i$ tel que
${f_{i1}}^2-4f_{i2}= (z_1^2-4z_2)R_i^2(z_1,z_2)$.   
Soit $\pi$ l'application de $\mathbb{C}^2$ dans
$\mathbb{C}^2$ avec $\pi(x,y)=(x+y,xy)$ ({\it voir} l'exemple 4). Posons
$h_i:=(h_{i1},h_{i2})=f_i\circ \pi$. On a
$h_{i1}^2-4h_{i2}=(x-y)^2R_i^2(x+y,xy)$. Ceci montre qu'il existe des
applications $F_i$ telles que $h_i=\pi\circ F_i$. Les $F_i$ sont
d\'efinies par:
$$F_i=(F_{i1},F_{i2})=\left(\frac{h_{i1}\pm (x-y)R_i(x+y,xy)}{2},
  \frac{h_{i2} \mp
(x-y)R_i(x+y,xy)}{2}\right).$$
On peut choisir les signes
 tels que $F_1\circ F_2=F_2\circ F_1$. Remarquons  
que $\tilde T_{d_i}(a_1+a_2)\circ \pi=a_1^{d_i}+a_2^{d_i}$. 
On v\'erifie facilement que
$h_{i1}=x^{d_i}+y^{d_i}+\o(|\nu|^{d_i})$ et
$h_{i2}=x^{d_i}y^{d_i}+\o(|\nu|^{d_i})$ o\`u $\nu:=(x,y)$. On a 
$F_{i1}=x^{d_i}+\o(|\nu|^{d_i})$ et
$F_{i2}=y^{d_i}+\o(|\nu|^{d_i})$. On a \'egalement
$F_{i1}(x,y)=F_{i2}(y,x)$. Utiliser le paragraphe pr\'ec\'edent
pour le couple $(F_1,F_2)$ on montre que ce couple est conjugu\'e \`a un
des couples des exemples 1-4. Comme aucune composante de 
$f_i{\cal C}_i$ n'est incluse dans ${\cal C}_i$, $F_i$
v\'erifie une propri\'et\'e analogue. De plus,
$F_i|_\infty=[x^{d_i}:y^{d_i}]$. Par cons\'equent,
il existe un automorphisme
$\sigma$ de
$\mathbb{P}^2$ qui pr\'eserve l'infini et qui v\'erifie
$\sigma^{-1}\circ F_i\circ\sigma=(\pm \tilde T_{d_i}x,\pm \tilde T_{d_i}y)$ ou 
$(\pm \tilde T_{d_i}y,\pm \tilde T_{d_i}x)$. Quitte \`a
remplacer $f_i$ par $f_i^4$, on peut supposer que  $\sigma^{-1}\circ
F_i\circ\sigma=(\tilde T_{d_i}x,\tilde T_{d_i}y)$. Comme
$F_i|_\infty=[x^{d_i}:y^{d_i}]$, on a
$\sigma=(\sigma_1,\sigma_2)=(\alpha x+\beta,\alpha' y+\beta')$. Comme
$F_{i1}(x,y)=F_{i2}(y,x)$, on a $\sigma_1^{-1}\circ \tilde
T_{d_i}\circ \sigma_1=\sigma_2^{-1}\circ\tilde T_{d_i}\circ \sigma_1$.
On a \'egalement 
$\sigma_1\{\pm
2\}=\sigma_2\{\pm 2\}$ car  l'ensemble de Julia
de $\tilde T_{d_i}$ est
$[-2,2]$. Ceci implique que $\beta=\beta'$ et
$\alpha=\pm\alpha'$. On montre facilement que $\alpha=-\alpha'$ est possible
seulement si les $d_i$ sont impairs et
$\beta=\beta'=0$. Dans ce cas, rien n'est chang\'e si l'on remplace
$\alpha'$ par $\alpha$. On peut supposer que
$\alpha=\alpha'$ et si $\beta=\beta'$. Il est clair maintenant que
$(f_1,f_2)$ est conjugu\'e \`a un couple construit dans
l'exemple 4.
\end{preuve}
\begin{proposition} La condition 2 est fausse.
\end{proposition}
\begin{preuve}  
Supposons que la condition 2 est vraie.
Consid\'erons d'abord le cas o\`u $d_1$ et $d_2$ sont impairs.
Alors
${\cal D}$ est la r\'eunion de ${\cal F}:=\{z_1=2\}$, de ${\cal
  F}'=\{z_1=-2\}$ et d'une courbe invariante $R$ de degr\'e 2 
passant par $\{b,c\}$ qui
est irr\'eductible ou une r\'eunion de deux droites. 
\par
Si $R$ est irr\'eductible, l'application $f_1$ agit sur $R$
comme un polyn\^ome de degr\'e $d_1$ et $b$, $c$ sont
exceptionnels pour $f_1|_R$. Par cons\'equent, $f_1|_R$ est
conjugu\'e \`a $x^{\pm d_1}$ et il
n'a pas de point critique diff\'erent de $b$ et $c$. Comme la
condition 2 est vraie, si $s\in \tilde T^{-1}\{\pm 2\}$ et
si $\{z_1=s\}$ coupe $R$ transversalement
en $z$, alors $z$ est un point critique de $f_1|_R$. C'est 
une contradiction.
\par
Alors $R$ est la r\'eunion de deux droites $L$ et $L'$
qui sont invariantes par $f_1$. La droite $L$ passe par $b$ et la
droite $L'$ passe par $c$. Quitte \`a effectuer un changement de
coordonn\'ees du type $(z_1,z_2)\mapsto (z_1,z_2+\alpha)$, on
peut suppser que $L$ est la droite $z_2-2z_1=0$. Il existe
$v\in\mathbb{C}$ tel que $L=\{z_2+2z_1+v=0\}$. D'apr\`es le lemme
6, il existe des polyn\^omes $T$ et $V$ tels que
\begin{equation}
f_{12}-2f_{11}=
(z_2-2z_1)T^2 \mbox{ et } f_{12}+2f_{11}+v=(z_2+2z_1+v)V^2
\end{equation}
On peut
\'ecrire 
$f_{12}=z_1^{d_1}\tilde T_{d_1}(z_2/z_1)+\Delta+\o(|z|^{d_1-1})$
o\`u $\Delta$ est homog\`ene de degr\'e $d_1-1$. On a
$f_{11}=\tilde T_{d_1}(z_1)=z_1^{d_1}+\o(|z|^{d_1-1})$ car il est
un polyn\^ome impair. Observons que les points critiques de
$\tilde T_{d_1}$ sont tous de multiplicit\'e 2 et s'envoient tous
dans les points fixes $\pm 2$. Alors il existe des polyn\^omes
homog\`enes $R$ et $S$ tels que $z_1^{d_1}\tilde
T_{d_1}(z_2/z_1)-2z_1^{d_1}= (z_2-2z_1)R^2$ et $z_1^{d_1}\tilde
T_{d_1}(z_2/z_1)+2z_1^{d_1}= (z_2+2z_1)S^2$. Il est clair que
$(z_2-2z_1)R$ et $(z_2+2z_1)S$ sont premiers entre eux. Alors on peut
\'ecrire $T=R+\Delta_1+\o(|z|^{(d_1-3)/2})$ et et
$V=S+\Delta_2+\o(|z|^{(d_1-3)/2})$ o\`u $\Delta_1$ et $\Delta_2$
sont homog\`enes de degr\'e $(d_1-3)/2$. D'apr\`es (4), on a
\begin{equation}
\Delta=2(z_2-2z_1)R\Delta_1=S[2(z_2+2z_1)\Delta_2+vS]
\end{equation} 
Ceci implique que $\Delta$ est divisible par $(z_2-2z_1)R$ et par
$S$. Comme $(z_2-2z_1)R$ et $S$ sont premiers entre eux, $\Delta$
est divisible par $(z_2-2z_1)RS$. Alors $\Delta=0$ car $\deg
\Delta<\deg (z_2-2z_1)RS$. On en d\'eduit que
$2(z_2+2z_1)\Delta_2+vS=0$. Alors $\Delta_2$ est divisible par
$S$. On obtient $\Delta_2=0$ et donc $v=0$ car $\deg \Delta_2<\deg S$.
Alors $L$ et $L'$ se coupent en $0$. Soit
$I:=(f_1|_L)^{-1}(0)\setminus\{0\}$.
Alors $f_1^{-1}L'$ passe par $I$. C'est la
contradiction recherch\'ee car $\#I=d_1-1$ et $\deg
f_1^{-1}L'=(d_1+1)/2$.   
\par
Supposons maintenant que $d_1$ est pair.
On peut supposer que $L=\{z_2-2z_1=0\}$ et $L'=\{z_2+2z_1+v=0\}$.
Il existe des polyn\^omes $R$ et $S$ tels que $z_1^{d_1}\tilde
T_{d_1}(z_2/z_1)-2z_1^{d_1}= (z_2-2z_1)(z_2+2z_1)R^2$ et
$z_1^{d_1}\tilde
T_{d_1}(z_2/z_1)+2z_1^{d_1}= S^2$. Il existe aussi des
polyn\^omes $T$ et $V$ tels que
$f_{12}-2f_{11}=(z_2-2z_1)(z_2+2z_1+v)T^2$ et
$f_{12}+2f_{11}+v=V^2$. Comme dans le cas pr\'ec\'edent, on
montre que
$v=0$ et que $f_1^{-1}L'$ passe par $I:=(f_1|_L)^{-1}(0)$. C'est
une contradiction car $\#I>\deg f_1^{-1}L'$.
\end{preuve}
\begin{proposition} Si la condition 1 est vraie, 
le couple $(f_1,f_2)$ est conjugu\'e \`a un couple
  d'applications polynomiales homog\`enes.
\end{proposition}
\begin{preuve} Quitte \`a remplacer $f_i$ par un de ses
it\'er\'es, on peut supposer que $d_i$ est suffisamment grand.
On consid\`ere le cas o\`u $d_1$, $d_2$ sont impairs. Les autres cas
seront trait\'es de m\^eme mani\`ere. 
D'apr\`es le lemme 6, on peut choisir un
  syst\`eme de coordonn\'ees tel que
${\cal D}=\{\Phi=0\}$ o\`u
$$\Phi=(z_2-2z_1)(z_2+2z_1)+u(z_2-2z_1)+v$$
est un polyn\^ome de degr\'e deux. 
Soit $0\leq \delta\leq d_1-1$ tel que on puisse \'ecrire
\begin{equation}
f_1=(f_{11},f_{12})=(z_1^{d_1},Q+\Delta+\o(|z|^\delta))
\end{equation}
o\`u $Q:=z_1^{d_1}\tilde 
T_{d_1}(z_2/z_1)$ et $\Delta$ est homog\`ene de degr\'e
$\delta$. D'apr\`es le lemme 6, il existe un polyn\^ome $W$ tel
que $\Phi(f_1)=\Phi. W^2$. Soient $R$ et $S$ des polyn\^omes
homog\`ens v\'erifiant $z_1^{d_1}\tilde
T_{d_1}(z_2/z_1)-2z_1^{d_1}=(z_2-2z_1)R^2$ et $z_1^{d_1}\tilde
T_{d_1}(z_2/z_1)+ 2z_1^{d_1}=(z_2-2z_1)S^2$. Alors $(z_2-2z_1)R$
et $(z_2+2z_1)S$ sont premiers entre eux.
\par
Montrons que $u=0$. On prend $\delta=d_1-1$. Quitte \`a remplacer
$W$ par $-W$ au cas n\'ecessaire, on peut \'ecrire
$W=RS+\Delta_W+ \o(|z|^{d_1-2})$ o\`u  $\Delta_W$ est homog\`ene
de degr\'e $d_1-2$. La relation $\Phi(f_1)=\Phi .W^2$ implique
\begin{eqnarray}
\lefteqn{[(z_2-2z_1)R^2+(z_2+2z_1)S^2]\Delta=} \nonumber \\
& = & RS[2(z_2-2z_1)(z_2+2z_1)\Delta_W
 +u(z_2-2z_1)RS]
\end{eqnarray} 
Alors $\Delta$ est divisible par $(z_2-2z_1)RS$. Comme $\deg
\Delta<(z_2-2z_1)RS$, on a $\Delta=0$ et donc
$2(z_2+2z_1)\Delta_W+uRS=0$. Alors $\Delta_W$ est divisible par
$RS$. On obtient $\Delta_W=0$ et donc $u=0$ car $\deg
\Delta_W<\deg RS$.
\par
De m\^eme mani\`ere, en utilisant $\delta=d_1-2$, on obtient
$v=0$.
\par
Soit $\delta$ l'entier minimal v\'erifiant (6). Alors on peut
\'ecrire 
$W=RS+\Delta_W+\o(|z|^{\delta-1})$ o\`u $\Delta_W$ est homog\`ene
de degr\'e $\delta-1$. On montre exactement comme
ci-dessus que $\Delta=0$. Comme $\delta$ est minimal, ceci
implique que $f_{12}$ est homog\`ene. 
\end{preuve}
\section{Preuve du th\'eor\`eme 1: troisi\`eme cas} 
Dans ce paragraphe, on suppose que les $f_i^n|_\infty$ ne sont pas
conjugu\'es ni \`a $z^{\pm  d_i^n}$ ni \`a
$\pm \tilde T_{d_1^n}$ ({\it voir} l'introduction). Il faut
montrer que $f_1$ et $f_2$ sont homog\`enes.
D'apr\`es \cite{Eremenko}, 
$f_i|_\infty$ d\'efinissent
des rev\^etements de $\cal{O}$ dans lui-m\^eme, 
o\`u ${\cal O}=(\mathbb{P}^1,n)$ est l'un
des orbifolds ${\cal O}_1,{\cal O}_2,{\cal O}_3,{\cal O}_4$ qui seront
d\'ecrits plus tard.  Il existe \'egalement un tore 
complexe $\mathbb{T}=\mathbb{C}/\Gamma$ de dimension 1, des
  applications lin\'eaires $\Lambda_i(x)=a_ix+b_i$
  pr\'eservant le groupe discret $\Gamma$ et une fonction
  elliptique $F:\mathbb{T}\rightarrow\mathbb{P}^1$ 
  tels que $f_i\circ F=F\circ \Lambda_i$. 
On a 
$d_i=|a_i|^2$ et donc tout point
p\'eriodique de
$f_i$ est r\'epulsif. En tout point fixe de $f_i$, la valeur propre
de $f_i$ est \'egale \`a $a_i$. 
L'orbifolds ${\cal O}_i$ est d\'efini par le couple
$(\mathbb{P}^1,n_i)$, o\`u $n_i$ est une application de $\mathbb{P}^1$
dans $\mathbb{Z}_+$
qui est \'egale \`a 1 sauf sur $s=3$ ou $s=4$ points. On
note $\alpha_j$ ces $s$ points. Alors l'ensemble des $\alpha_j$
est invariant par $f_1$ et $f_2$. Quitte \`a remplacer $f_1$,
$f_2$ par $f_1\circ f_2^n$ et $f_1\circ f_2^m$ avec
$n\not=m$ convenables, on peut supposer que les
$\alpha_j$ qui sont p\'eriodiques pour $f_1$, sont \'egalement
p\'eriodiques
pour $f_2$ et r\'eciproquement. En effet, si $f_1\circ f_2^n$ et
$f_1\circ f_2^m$ sont homog\`enes, $f_1$ et $f_2$ le sont aussi.
Maintenant, quitte \`a remplacer
$f_i$ par l'un de ses it\'er\'es, on peut supposer que si
$\alpha_j$ est p\'eriodique, il est fixe.   
Par d\'efinition des orbifolds et
leur description \cite{Eremenko},
on a les cas suivants:
\par {\bf 1.} Pour ${\cal O}_1$, $s=3$ et 
$n_1(\alpha_j)=3$ pour $j=1$, 2, 3 et
$a_i\in\mathbb{Z}[\omega]$ avec $\omega=\exp(\pi i/3)$. Alors
si $d_i$ n'est pas
  divisible par 3, $d_i=1$ modulo 3. 
  Les $\alpha_j$ sont tous fixes r\'epulsifs 
  pour $f_i|_\infty$
  et l'image r\'eciproque de $\alpha_j$ 
 est constitu\'e par lui-m\^eme et $(d_i-1)/3$ autres 
 points; la multiplicit\'e  de $f_i|_\infty$ sur chacun de ces
$(d_i-1)/3$ 
 points est \'egale \`a 3. 
\par
Sinon,
  $f_i(\alpha_1)=f_i(\alpha_2)=f_i(\alpha_3)=\alpha_1$. 
L'image r\'eciproque de $\alpha_j$
  pour $j=2,3$ est constitu\'e par $d_i/3$ points  de
  multiplicit\'e 3;  celle de $\alpha_1$
  est constitu\'e par les $\alpha_j$ et $d_i/3-1$ autres points  de
  multiplicit\'e 3.
\par
{\bf 2.} Pour ${\cal O}_2$, $s=3$ et $n_2(\alpha_1)=6$,
$n_2(\alpha_2)=3$, 
$n_2(\alpha_3)=2$ et
$a_i\in\mathbb{Z}[\omega]$ avec $\omega=\exp(\pi i/3)$.
Alors si $d_i$ n'est pas
  divisible par 2 ni par 3, $d_i=1$ modulo 6. 
  Les $\alpha_j$ sont tous fixes r\'epulsifs 
  pour $f_i$
  et l'image r\'eciproque de $\alpha_1$ (resp. $\alpha_2$ et $\alpha_3$) 
 est constitu\'e par lui-m\^eme et $(d_i-1)/6$ (resp.  $(d_i-1)/3$ 
et  $(d_i-1)/2$) autres
 points  de multiplicit\'e 6 (resp. 3 et 2). 
\par
Si $d_i$ est 
  divisible par 2 mais non pas par 3, alors $d_i=4$ modulo 6. Dans ce cas,
  $\alpha_1$ et $\alpha_2$ sont fixes et $f_i(\alpha_3)=\alpha_1$.  
L'image r\'eciproque de $\alpha_3$
  est constitu\'e par $d_i/2$ points  de
  multiplicit\'e 2;  celle de $\alpha_2$ (resp. $\alpha_1$)
  est constitu\'e par lui-m\^eme (resp. lui-m\^eme et $\alpha_3$ avec
  la multiplicit\'e 3) 
et $(d_i-1)/3$ (resp. $(d_i-4)/6$) autres points  de
  multiplicit\'e 3 (resp. 6).
\par
Si $d_i$ est 
  divisible par 3 mais non pas par 2, alors $d_i=3$ modulo 6. Dans ce cas,
  $\alpha_1$ et $\alpha_3$ sont fixes et $f_i(\alpha_2)=\alpha_1$.  
L'image r\'eciproque de $\alpha_2$
  est constitu\'e par $d_i/3$ points  de
  multiplicit\'e 3;  celle de $\alpha_3$ (resp. $\alpha_1$)
  est constitu\'e par lui-m\^eme (resp. lui-m\^eme et $\alpha_2$ avec
  la multiplicit\'e 2) 
et $(d_i-1)/2$ (resp. $(d_i-3)/6$) autres  points  de
  multiplicit\'e 2 (resp. 6).
\par
 Si $d_i$ est divisible par 6,
  $f_i(\alpha_1)=f_i(\alpha_2)=f_i(\alpha_3)=\alpha_1$. L'image
  r\'eciproque de $\alpha_3$ (resp. $\alpha_2$) est constitu\'e par
  $d_i/2$ (resp. $d_i/3$) points  de multiplicit\'e 2
  (resp. 3); celle de $\alpha_1$  est constitu\'e par lui-m\^eme,
  $\alpha_3$ avec la multiplicit\'e 3, $\alpha_2$ avec la
  multiplicit\'e 2 et
  $d_i/6$ points  de multiplicit\'e 6.
\par
{\bf 3.} Pour ${\cal O}_3$, $s=3$ et $n_1(\alpha_1)=4$,
$n_1(\alpha_2)=4$, $n_1(\alpha_3)=2$
et $a_i\in\mathbb{Z}[i]$. Alors si $d_i$ n'est pas
  divisible par 2, $d_i=1$ modulo 4. Dans ce
  cas, les
  $\alpha_j$ sont fixes r\'epulsifs.
L'image r\'eciproque de $\alpha_3$ (resp. $\alpha_2$ ou $\alpha_1$) 
 est constitu\'e par lui-m\^eme et $(d_i-1)/2$ (resp. $(d_i-1)/4$) 
 autres points  de multiplicit\'e 2 (resp. 4). 
\par
Sinon, quitte \`a remplacer $f_i$ par $f_i^2$, on peut supposer
que $d_i$ est divisible par 4. Alors
  $f_i(\alpha_1)=f_i(\alpha_2)=f_i(\alpha_3)=\alpha_1$. 
L'image r\'eciproque de $\alpha_3$ (resp. $\alpha_2$)
  est constitu\'e par $d_i/2$ (resp. $d_i/4$) points  de
  multiplicit\'e 2 (resp. 4);  celle de $\alpha_1$
  est constitu\'e par lui-m\^eme, $\alpha_2$ avec multiplicit\'e 1,
  $\alpha_3$ avec multiplicit\'e 2 et $d_i/4-1$ autres points  de
  multiplicit\'e 4.
\par
{\bf 4.} Pour ${\cal O}_4$, $s=4$ et $n_1(\alpha_j)=2$ pour $j=1,2,3,4$. 
Si $d_i$ n'est pas
  divisible par 2, alors tous les 
$\alpha_j$ sont fixe r\'epulsifs. 
L'image r\'eciproque de chaque $\alpha_j$ 
 est constitu\'e par lui-m\^eme et $(d_i-1)/2$ autres  
 points  de multiplicit\'e 2. 
\par
Sinon, il y a deux cas possibles.
Pour le premier cas,
  $f_i(\alpha_1)=f_i(\alpha_2)=f_i(\alpha_3)=f_i(\alpha_4)=\alpha_1$. 
L'image r\'eciproque de $\alpha_j$ pour $j=2,3,4$
  est constitu\'e par $d_i/2$ points  de
  multiplicit\'e 2; celle de $\alpha_1$
  est constitu\'e par les $\alpha_j$ avec la multiplicit\'e 1 
  et $d_i/2-2$ autres points  de
  multiplicit\'e 2.
\par
Pour le second,
  $f_i(\alpha_1)=f_i(\alpha_3)=\alpha_1$ et
$f_i(\alpha_2)=f_i(\alpha_4)=\alpha_2$. 
L'image r\'eciproque de $\alpha_j$ pour $j=3,4$
  est constitu\'e par $d_i/2$ points  de
  multiplicit\'e 2; celle de $\alpha_j$ pour $j=1,2$ 
  est constitu\'e par $\alpha_j$, $\alpha_{j+2}$
avec la multiplicit\'e 1  et $d_i/2-1$ autres points  de
  multiplicit\'e 2.
\par
Remarquons que dans tous les cas, on a
$$\sum_j\frac{n(\alpha_j)-1}{n}=2.$$
Soit $r$
le nombre des $\alpha_j$ qui sont fixes par $f_1$ et $f_2$.
Alors $r=1$, $2$ ou $s$. 
\begin{lemme} Si $r=s$, il existe une courbe
 $L$ de $\mathbb{P}^2$ invariante par $f_1$ et $f_2$; elle
  coupe $\{w_0=0\}$ transversalement en
  $\alpha_1,\ldots,\alpha_s$.
De plus,  
$f_i^{-1}L=L\cup {\cal C}_i$ et $\mult(f_i,A)$  
  \'egale \`a la multiplicit\'e de 
$f_i|_\infty$ en chaque point de $A\cap\{w_0=0\}$ 
pour toute composante $A$ de ${\cal C}_i$.
\end{lemme}
\begin{preuve} Il suffit de montrer la proposition pour $i=1$. 
Soit $U$ un voisinage assez petit de $\{w_0=0\}$ tel
  que ${\cal C}_1\cap U$ soit la r\'eunion des courbes
  irr\'eductibles $A_p$; chaque $A_p$ coupe $\{w_0=0\}$ en un seul point
  $\beta_p$. D'apr\`es le lemme 3, $\beta_p$ appartient \`a l'ensemble
  critique de $f_1|_\infty$. Par chaque point $\alpha_j$ passe une 
  courbe lisse, stable par $f_1$ et par $f_2$. Fixons un $A_p$. Soit
  $\alpha_j=f_1(\beta_p)$. On choisit 
un syst\`eme de coordonn\'ees
  local $(x,y)$ pour un voisinage de $\alpha_j$ tel que
  $\alpha_j=(0,0)$, $\{x=0\}=\{w_0=0\}$ soit 
  la vari\'et\'e instable, $(y=0)$ soit
  la vari\'et\'e stable et tel que
  $f_i(x,y)=(x^{d_i}+\o(x^{d_i}),a_iy+xyg_i(x,y))$. Comme
$(f_1|_\infty)^n\not=(f_2|_\infty)^m$
pour tous entiers positifs $n$, $m$,
  on a $a_1^na_2^m\not = a_1^{n'}a_2^{m'}$ pour tous
  $(n,m)\not=(n',m')$. On va montrer que
  $f_1A_p$ est la vari\'et\'e stable de $f_1$ en $\alpha_j$.
Il est clair que $f_1A_p\not =\{x=0\}$. D'apr\`es la proposition 1,
  $f_1^nf_2^mA_p=f_1^{n'}f_2^{m'}A_p$ pour certains couples $(n,m)\not
  =(n',m')$. D'o\`u $f_1^nf_2^m(f_1A_p)=f_1^{n'}f_2^{m'}(f_1A_p)$. 
D'apr\`es le lemme 5, $f_1A_p$ est la vari\'et\'e
  stable car elle passe par $\alpha_j$. 
\par
Soit $L:=f_1({\cal C}_1)$.
Comme les vari\'et\'es stables coupent $\{w_0=0\}$ transversalement,
  $L$ coupe $\{w_0=0\}$ transversalement en les $\alpha_j$. Donc $\deg
  L=s$.  Alors les $A_p$ sont les composantes
  irr\'eductibles de $\overline{f_1^{-1}L\setminus L}\cap U$. Soient
  $\beta_p:=A_p\cap\{w_0=0\}$ et $\alpha_{j_p}:=f_1(\beta_p)$ avec
  $1\leq j_p\leq s$. 
Posons $m_p:=\mult(f_1,A_p)$ et 
  $n_p:=\mult(A_p\cap\{w_0=0\},\beta_p)$. D'apr\`es le lemme 4,
  $\sum_{j_p=j}m_pn_p=d_1-1$,
$\sum m_pn_p=s(d_1-1)$ et $m_p\leq
\mult(f_1|_\infty,\beta_p)=n(\alpha_{j_p})$.
On sait que $\deg {\cal C}_1=2(d_1-1)$ compt\'e avec les
multiplicit\'es. D'apr\`es le lemme 3,  $\sum
(m_p-1)n_p=2(d_1-1)$. D'o\`u:
\begin{eqnarray*}
\sum_p(m_p-1)n_p & = &\sum_p\frac{m_p-1}{m_p}m_pn_p \leq  \sum_p
\frac{n(\alpha_{j_p})-1}{n(\alpha_{j_p})}m_pn_p \\
& = & 
\sum_j\frac{n(\alpha_j)-1}{n(\alpha_j)}\sum_{j_p=j}m_pn_p
 = \sum_j\frac{n(\alpha_j)-1}{n(\alpha_j)}(d_1-1)\\
&  = &  2(d_1-1)
\end{eqnarray*}
Par cons\'equent, $m_p=n(\alpha_{j_p})=\mult(f_1|_\infty,\beta_{p})$
et $n_p=1$ 
pour tout $p$. Alors les $A_p$ sont disjointes. 
Alors
il est clair que ${\cal C}_1=\overline{f_1^{-1}L\setminus L}$ et
$\mult(f_1,A)$ est \'egale \`a la multiplicit\'e de $f_1|_\infty$ en
chaque point de $A\cap\{w_0=0\}$ pour toute composante $A$ de $f_1^{-1}L$.
\end{preuve}
\begin{lemme} Si $r=1$, il existe une droite  $L\not
  =\{w_0=0\}$ invariante par
  $f_1$, $f_2$ et passant par $\alpha_1$. L'ensemble $f_i^{-1}L$ se
  constitue par $L$, par des courbes $L'$ et 
  $H_i\subset {\cal C}_i$. 
La courbe $L'$ coupe $\{w_0=0\}$ transversalement en
  $\alpha_2,\ldots, \alpha_s$ et $f_i^{-1}L'=\overline{{\cal
  C}_i\setminus (H\cup L')}$. Pour chaque composante $A$
  de ${\cal C}_i$, $\mult(f_i,A)$ est \'egale \`a la multiplicit\'e de
  $f_i|_\infty$ en chaque point de $A\cap\{w_0=0\}$.
\end{lemme}
\begin{preuve} On peut supposer que $i=1$.
On consid\`ere le cas o\`u ${\cal O}={\cal O}_2$. Les autres cas
seront trait\'es de m\^eme mani\`ere.
\par
Soient $U$, $A_p$, $\beta_p$ d\'efinies dans la preuve du
lemme pr\'ec\'edent. Fixons un $p$ tel que $f_1A_p=\alpha_1$. 
Comme dans le lemme pr\'ec\'edent, on montre que
$f_1A_p$ est la vari\'et\'e stable de $f_1$ en $\alpha_1$. Soit $L$
la courbe alg\'ebrique contenant $f_1A_p$. Alors $L$ coupe $\{w_0=0\}$
transversalement en un point unique $\alpha_1$ car ${\cal
C}_1\cap \{w_0=0\}$ s'envoie dans $\alpha_1$.
La courbe $L$ est donc une
droite. Soient $B_q$, $C_\nu$, $D_\xi$ les composantes irr\'eductibles de
$\overline{f_1^{-1}L\setminus L}\cap U$ telles que
$C_\nu\cap\{w_0=0\}=\alpha_2$, $D_\xi\cap\{w_0=0\}=\alpha_3$  et 
$B_q\cap\{w_0=0\}=:\beta_q\not=\alpha_2,\alpha_3$. On appele $m_q$, $m_\nu$,
$m_\xi$ les multiplicit\'es de $f_1$ sur $B_q$, $C_\nu$, $D_\xi$ et 
$n_q$, $n_\nu$,
$n_\xi$ les multiplicit\'es de l'intersection de  $B_q$, $C_\nu$,
$D_\xi$ avec $\{w_0=0\}$.
\par
D'apr\`es les lemmes 4 et 3, on a 
$\sum m_qn_q=d_1-6$, $\sum m_\nu n_\nu=2$, $\sum m_\xi n_\xi=3$, 
$\sum (m_q-1)n_q=5(d_1/6-1)$, $\sum (m_\nu-1) n_\nu=1$, 
$\sum (m_\xi-1) n_\xi=2$, $m_q\leq 6$, $m_\nu\leq 2$ et
$m_\xi\leq 3$. Donc
$$5\left(\frac{d_1}{6}-1\right)=\sum(m_q-1)n_q =
\sum\frac{m_q-1}{m_q}m_qn_q\leq\sum\frac{5}{6}m_qn_q = 
5\left(\frac{d_1}{6}-1\right).$$
D'o\`u $m_q=6=\mult(f_1|_\infty,\beta_q)$ et $n_q=1$. 
De m\^eme, $m_\nu=2=\mult(f_1|_\infty,\alpha_2)$, $n_\nu=1$,
$m_\xi=3=\mult(f_1|_\infty,\alpha_3)$ et 
$n_\xi=1$. Il y a donc $d_1/6-1$ composantes $B_q$;
sur chacune $B_q$ la multiplicit\'e de $f_1$ est \'egale \`a 6; 
il y a une seule composante $C_\nu$ 
avec la multiplicit\'e 2 et une seule composante $D_\xi$ avec la
multiplicit\'e 3. Ces composantes sont disjointes et coupent
$\{w_0=0\}$ transversalement.
Posons $C:=C_\nu$ et $D:=D_\xi$. Soit $A_p$ une composante v\'erifiant
$f_1A_p\cap \{w_0=0\}=\alpha_2$. Comme dans
le lemme 8, on prouve que $f_1^2A_p$ est la vari\'et\'e stable
de $f_1$ en $\alpha_1$. D'o\`u $f_1A_p=C$. Soient $F_\eta$ les
composantes irr\'eductibles de $f_1^{-1}C\cap U$, $\beta_\eta:=
F_\eta\cap\{w_0=0\}$, $m_\eta:=\mult(f_1,F_\eta)$ et $n_\eta:=
\mult(F_\eta\cap\{w_0=0\},\beta_\eta)$. D'apr\`es les lemmes 3 et
4,  on a $\sum m_\eta
n_\eta=d_1$, $\sum (m_\eta-1)n_\eta=2d_1/3$ et $m_\eta \leq
3$. D'o\`u:
$$\frac{2}{3}d_1=\sum (m_\eta-1)n_\eta=\sum\frac{m_\eta-1}{m_\eta}m_\eta
n_\eta\leq\sum\frac{2}{3} m_\eta n_\eta=\frac{2}{3}d_1.$$
Ceci implique que $n_\eta=1$,
$m_\eta=3=\mult(f_1|_\infty,\beta_\eta)$. 
Les $F_\eta$ appartiennent
\`a ${\cal C}_1$; chacune est de multiplicit\'e 3 et coupe
$\{w_0=0\}$ transversalement.
La courbe $f_1^{-1}D$ v\'erifie une propri\'et\'e analogue.
Alors $f_1^{-1}L'=\overline{{\cal C}_1\setminus(H\cup L')}$ et pour
chaque composante $A$ de ${\cal C}_1$, $\mult(f_1,A)$ est \'egale \`a
la multiplicit\'e de $f_1|_\infty$ en chaque point de $A\cap \{w_0=0\}$.
\end{preuve}
\begin{lemme}
Si $r=2$ et si ${\cal O}={\cal O}_2$,
il existe deux droites $L$ et $L'$ invariantes par $f_1$ et $f_2$;
$L$ passe par $\alpha_1$ et $L'$ passe par $\alpha_j$ avec $j=2$ ou
$3$.
La courbe $f_i^{-1}L$ est la r\'eunion de $L$, d'une droite
$\Lambda$ passant par $\alpha_l$ et d'une autre courbe $H_i$ o\`u
$l\not=1,j$.
La courbe $f_i^{-1}L'$ est la r\'eunion de $L'$ 
et d'une autre courbe $H_i'$.
On a ${\cal C}_i=\Lambda\cup f_i^{-1}\Lambda \cup H_i \cup H_i'$.
Pour toute composante
irr\'eductible $A$ de ${\cal C}_i$, $\mult(f_i,A)$ est \'egale \`a la
multiplicit\'e de $f_i|_\infty$ en chaque point de $A\cap\{w_0=0\}$.
\end{lemme}
\begin{preuve} 
Si $d_i$ est pair, $j=2$ et
  $l=3$; si $d_i$ est divisible par 3, $j=3$ et $l=2$.
Comme dans les lemmes pr\'ec\'edents, on montre qu'il existe
  une courbe invariante passant par $\alpha_1$ et $\alpha_j$. L'image
  r\'eciproque de cette courbe par $f_i$ contient des composantes de
  ${\cal C}_i$, o\`u les multiplicit\'es de $f_i$ sont \'egales \`a 
  $n(\alpha_1)$ et $n(\alpha_j)$. Comme $n(\alpha_1)\not=n(\alpha_j)$,
  cette courbe
  invariante est r\'eductible. Elle est donc la r\'eunion de deux
  droites $L$ et
  $L'$. On montre tout comme dans les lemmes pr\'ec\'edents que
  $f_i^{-1}L$ contient $\Lambda$ passant par $\alpha_l$, que ${\cal
  C}_i=\Lambda\cup f_i^{-1}\Lambda\cup H_i \cup H_i'$
et que
  $\mult(f_i,A)$ est \'egale \`a la multiplicit\'e de
$f_i|_\infty$ en chaque point de
$A\cap\{w_0=0\}$ pour toute composante $A$ de ${\cal
  C}_i$.
\end{preuve}
\begin{lemme}
Si $r=2$ et si ${\cal O}_4$,
il existe une courbe ${\cal D}$
invariante par $f_1$ et $f_2$; cette courbe est de degr\'e 2 et
passe par $\alpha_1$ et $\alpha_2$.
La courbe $f_i^{-1}{\cal D}$ est la r\'eunion de ${\cal D}$,
d'une courbe
${\cal D}'$ de degr\'e 2 passant par $\alpha_3$ et $\alpha_4$
et d'une autre courbe $H_i$.
On a ${\cal C}_i=H_i \cup f_i^{-1}{\cal D}'$.
Pour toute composante
irr\'eductible $A$ de ${\cal C}_i$, $\mult(f_i,A)$ est \'egale \`a 2.
\end{lemme}
\begin{preuve} 
Comme dans
les lemmes pr\'ec\'edents, on montre qu'il existe une courbe
${\cal D}$ invariante par $f_1$ et $f_2$ et coupe $\{w_0=0\}$
transversalement en $\alpha_1$ et $\alpha_2$. Cette courbe est
alors de degr\'e 2.
L'ensemble
$f_i^{-1}({\cal D})$ est la r\'eunion d'une courbe ${\cal D}'$ et
d'une courbe ${\cal H}_i$. La courbe ${\cal D'}$ coupe $\{w_0=0\}$
transversalement en $\alpha_3$ et $\alpha_4$. La courbe ${\cal
D}'$ est \'egalement de degr\'e 2.
On a aussi ${\cal C}_i={\cal H}_i\cup f_i^{-1}{\cal D}'$.
On  montre comme dans les lemmes pr\'ec\'edents que   
  $\mult(f_i,A)$ est \'egale \`a la multiplicit\'e de $f_i|_\infty$
en chaque point de $A\cap\{w_0=0\}$ (i.e. \'egale \`a 2)
pour toute composante $A$ de ${\cal  C}_i$.
\end{preuve}
\begin{proposition} Le couple
  $(f_1,f_2)$ est conjugu\'e \`a un couple d'applications
 polynomiales homog\`enes.
\end{proposition}
\begin{preuve} Nous allons traiter seulement le cas o\`u $r=s$ et
${\cal O}={\cal O}_4$ qui nous semble le plus compliqu\'e.
La preuve est essentiellement valable pour les autres cas.
\par
Quitte \`a remplacer $f_1$ par l'un de ses it\'er\'es, on peut
supposer que $d_1$ est suffisamment grand.
On choisit un syst\`eme de coordonn\'ees tel que
  $\alpha_1=[0:1:0]$, $\alpha_2=[0:0:1]$, $\alpha_3=[0:1:1]$,
$\alpha_4=[0:1:\alpha]$
et la courbe $L$ soit d\'efinie par l'\'equation $\Phi=0$ o\`u
$\alpha\in \mathbb{C}\setminus\{0,1\}$ et
$$\Phi(z)=
z_1z_2(z_2-z_1)(z_2-\alpha z_1)+bz_1^2z_2+cz_1z_2^2+dz_1+uz_2+v$$
est un polyn\^ome de degr\'e 4.
Soit $0\leq\delta\leq d_1-1$ un nombre naturel tel qu'on peut
\'ecrire 
\begin{equation}
f_1(z)=(f_{11},f_{12})=
(P+\Delta_1+\o(|z|^{\delta}),Q+\Delta_2+\o(|z|^{\delta}))
\end{equation} 
o\`u  $P,Q$ (resp. $\Delta_1$ et $\Delta_2$)
  sont homog\`enes de degr\'e $d_1$ (resp. $\delta$).
D'apr\`es le lemme 8, il existe un polyn\^ome $W$ tel que
$\Phi(f_1)=\Phi.W^2$. On pose $\Delta_3:=\Delta_2-\Delta_1$ et
$\Delta_4:=\alpha\Delta_2 -\Delta_1$. Par d\'efinition de ${\cal
O}_4$, il existe des polyn\^omes homog\`enes $R$, $S$, $T$ et $V$
de degr\'e $(d_1-1)/2$
tels que $P=z_1R^2$, $Q=z_2S^2$, $Q-P=(z_2-z_1)T^2$ et $\alpha
Q-P=(\alpha z_2-z_1)V^2$. Observons que $z_1R$, $z_2S$,
$(z_2-z_1)T$ et $(\alpha z_2-z_1)V$ sont deux \`a deux premiers
entres eux car $P$ et $Q$ sont permiers entre eux.
On a  $W=RSTV+\o(|z|^{2(d_1-1)})$.
\par
Montrons que $b=c=0$. Consid\'erons $\delta=d_1-1$. On \'ecrit
$W=RSTV+ \Delta+\o(|z|^{2d_1-3})$ o\`u $\Delta$ est homog\`ene de
degr\'e $2d_1-3$. On obtient de l'\'equation
$\Phi(f_1)= \Phi W^2$ que
\begin{eqnarray}
\lefteqn{R^2S^2T^2\Delta_4+R^2S^2V^2\Delta_3+
R^2T^2V^2\Delta_2+S^2T^2V^3\Delta_1} &   & \nonumber\\
& = &
RSTV[(bz_1^2z_2+cz_1z_2^2) RSTV +2 z_1 z_2
(z_2-z_1)(z_2-\alpha z_1)\Delta]
\end{eqnarray}
Cette \'egalit\'e implique que $\Delta_1$ (resp. $\Delta_2$,
$\Delta_3$ et $\Delta_4$) est divisible par $z_1R$ (resp. $z_2S$,
 $T$ et $V$). On pose $\Delta_1=z_1R\Delta_1'$,
$\Delta_2=z_2S\Delta_2'$, $\Delta_3=T\Delta_3'$ et
$\Delta_4=V\Delta_4'$. Par d\'efinition de $\Delta_3$ et de $T$
on a
\begin{eqnarray}
z_2S\Delta_2'-z_1R\Delta_1' & = & T\Delta_3'\\
z_2S^2-z_1R^2 & = & (z_2-z_1)T^2
\end{eqnarray}
On obtient de ces relations que
\begin{equation}
z_1R(S\Delta_1'-R\Delta_2')=T[(z_2-z_1)T\Delta_2'-S\Delta_3']
\end{equation}
Par cons\'equent, $S\Delta_1'-R\Delta_2'$ est divisible par $T$
car $z_1R$ et $T$ sont premiers entre eux.
De m\^eme mani\`ere, \`a l'aide de $\Delta_4$ et de $V$, on
montre que $S\Delta_1'-R\Delta_2'$ est divisible par $V$. Comme
$T$ et $V$ sont premiers entre eux, $S\Delta_1'-R\Delta_2'$ est
divisible par $TV$. Mais $\deg S\Delta_1'-R\Delta_2'$ est plus
petit que $\deg TV$. Donc $S\Delta_1'-R\Delta_2'=0$. Ceci
implique que  $\Delta_1'$ est divisible par $R$. On en d\'eduit
que $\Delta_1'=\Delta_2'=0$ car $\deg \Delta_1'<\deg R$. On
obtient aussi $\Delta_1=\Delta_2=0$ et $b=c=0$.
\par
De m\^eme mani\`ere, avec $\delta=d_1-2$, on montre que $d=u=0$ et
ensuite avec $\delta=d_1-3$, on montre que $v=0$.
\par  
On choisit $\delta$ minimal tel que (8) soit vrai. Alors
$\delta\leq d_1-4$. 
D'apr\`es l'\'equation $\Phi(f_1)=\Phi
W^2$, on peut \'ecrir $W=RSTV+\Delta+\o(|z|^{d_1-2+\delta})$
o\`u $\Delta$ est homog\`ene de degr\'e $d_1-2+\delta$. On montre
exactement comme ci-dessus que $\Delta_1=\Delta_2=0$. Comme
$\delta$ est minimal, ceci implique que $f_1$ est homog\`ene. De
m\^eme mani\`ere, on montre que $f_2$ est homog\`ene.
\end{preuve}

\begin{thebibliography}{11}
\bibitem{Eremenko}
\textit{A.E. Eremenko}, On some functional equations connected with
iteration of rational function, \textit{Leningrad. Math. J.}, \textbf{1}
(1990), No. 4, 905-919.
%
\bibitem{Fatou} 
\textit{P. Fatou}, Sur l'it\'eration analytique et les substitutions
permutables, \textit{J. Math.}, \textbf{2} (1923), 343. 
%
\bibitem{FornessSibony}
\textit{J.E. Fornæss, N. Sibony}, Complex dynamics in higher
dimension I, \textit{Ast\'erique}, \textbf{222} (1994), 201-213.
% 
%
%
\bibitem{Julia}
\textit{G. Julia}, M\'emoire sur la permutabilit\'e des fractions
rationnelles, \textit{Ann. Sci. Ecole Norm. Sup.}, \textbf{39} (1922),
131-215.
%
\bibitem{Lamy}
\textit{S. Lamy}, Alternative de Tits pour ${\rm Aut}[\mathbb{C}^2]$,
{\it Pr\'epublication}. 
%
%
\bibitem{LevinPrzytycki}
\textit{G. Levin, F. Przytycki}, When do two functions have the same
Julia set?, \textit{Proc. Amer. Math. Soc.} \textbf{125} (1997),
no. 7, 2179-2190. 
%
\bibitem{Ritt}
\textit{J.F. Ritt}, Permutable rational functions,
\textit{Trans. Amer. Math. Soc}, \textbf{25} (1923), 399-448.
%
%

\bibitem{Sibony}
\textit{N. Sibony}, Dynamique des applications rationelles de
$\mathbb{P}^k$, \textit{Survey}, (1999).
%
\bibitem{Veselov}
\textit{A.P. Veselov}, Integrable mappings and Lie algebras,
\textit{Dokl. Akad. Nauk SSSR} \textbf{292} (1987), 1289-1291; English
transl. in Soviet Math. Dokl., \textbf{35} (1987).
\end{thebibliography}
\end{document}